\documentclass[11pt]{amsart}
\usepackage{amssymb}
\usepackage{color}
\usepackage{verbatim} 
\newtheorem{Theorem}[equation]{Theorem}
\newtheorem{Proposition}[equation]{Proposition}
\newtheorem{cor}[equation]{Corollary}
\newtheorem{Lemma}[equation]{Lemma}

\newtheorem{Definition}[equation]{Definition}
\newtheorem{rem}[equation]{Remark}

\numberwithin{equation}{section}

\address{%
}

\def\a1s{a_1,\cdots, a_s}
\def\a{\alpha}

\def\aa{\mathcal A}

\def\fa{\mathfrak{a}}

\def\andd{\quad\hbox{and}\quad}

\def\fb{\frak{b}}
\def\b{\beta}

\def\bb{\mathcal{B}}

\def\bl4{B_{\ell\geq4}}
\def\cc{{\mathcal C}}

\def\d{\delta}
\def\D{\Delta}

\def\dd{\mathcal D}

\def\ve{\varepsilon}

\def\bbbf{\mathbb{F}}

\def\gg{{\mathcal G}}

\def\fg{\mathfrak{g}}

\def\heart{\hbox{\tiny$\heartsuit$}}

\def\hh{{\mathcal H}}

\def\fh{\mathfrak{h}}

\def\ii{\mathcal{I}}
\def\jj{\mathcal{J}}

\def\kk{\mathcal{K}}

\def\lam{\lambda}
\def\Lam{\Lambda}
\def\LL{\mathcal{L}}

\def\ep{\epsilon}
\def\fm{(\cdot,\cdot)}

\def\m{\mathcal{M}}
\def\bbbn{\mathbb{N}}

\def\1k{\frac{1}{k}}
\def\op{\oplus}
\def\ot{\otimes}

\def\la{\langle}
\def\ra{\rangle}

\def\sub{\subseteq}
\def\sg{\sigma}

\def\rcross{R^{\times}}

\def\pf{\noindent{\bf Proof. }}

\def\i{{\mathcal I}}

\def\kk{\mathcal K}

\def\sfs{\mathfrak{s}}

\def\ss{\mathcal{S}}

\def\T{{\mathcal T}}

\def\u{{\mathcal U}}

\def\v{{\mathcal V}}

\def\w{{\mathcal W}}

\def\bbbz{{\mathbb Z}}
\def\ss{\mathcal{S}}
 
\def\1il{1\leq i\leq\ell}

\begin{document}

\centerline{\bf Structure of root graded Lie algebras}

\vspace{1cm}
\centerline{Malihe Yousofzadeh\footnote{This research was in part
supported by a grant from IPM (No. 89170030).}}


\vspace{1cm}
\parbox{5in}{Abstract. We give a complete description of Lie algebras graded by an infinite irreducible locally finite root system.}
%

\section{introduction}
In 1992, S. Berman and R. Moody \cite{BM} introduced the notion of
a  Lie algebra graded by an irreducible  reduced finite  root
system. Their   definition was  motivated by a construction
appearing  in the classification of finite dimensional simple Lie
algebras containing  nonzero  toral subalgebras \cite{Se}. The
classification of root graded Lie algebras in the sense of S.
Berman and R. Moody was given,  in part,  by S. Bermen and R.
Moody themselves and was completed   by G. Benkart and E. Zelmanov
\cite{BZ} in 1996. This classification has  been based on a
type-by-type approach; for each type $X,$ the authors give a
recognition theorem for centerless Lie algebras graded by a root
system of type $X.$  In 1996, E. Neher \cite{N} generalized the
notion  of root graded Lie algebras by switching from fields of
characteristic zero to  rings containing $1/6$ and working with
locally finite root systems instead of  finite root systems.
Roughly speaking, according to him,  a Lie algebra $\LL$ over a ring
containing $1/6$  is  graded by a reduced locally finite root
system $R$ if $\LL$ is a $Q(R)-$graded Lie algebra generated by
homogenous submodules of nonzero degrees and that for any nonzero
root $\a\in R,$ there are homogenous elements $e$ and $f$ of
degrees $\a$ and $-\a$ respectively such that  $[e,f]$ acts
diagonally on $\LL.$  He realized root graded Lie algebras for
reduced  types other than $F_4,$ $G_2$ and $E_8$ as  central
extensions  of  Tits-Kantor-Koecher  algebras of certain Jordan
pairs. Finally in 2002, B. Allison, G. Benkart and Y. Gao
\cite{ABG2}  defined a Lie algebra   graded by an irreducible
finite root system of type $BC$  and studied  root graded Lie
algebras of type $BC_n$ for $n\geq2.$  In 2003, G. Benkart and O.
Smirnov \cite{BS}  studied Lie algebras graded by a finite root
system of type $BC_1$ and  finalized the classification of Lie
algebras graded by an irreducible  finite root system.

A  Lie algebra $\LL$ graded by an irreducible  finite root system
$R$ has a weight space decomposition with respect to a splitting
Cartan subalgebra of a finite dimensional split simple Lie
subalgebra $\fg$ of $\LL,$  whose  set of weights is contained  in
$R.$ This feature allows us to decompose $\LL$ as $\LL=\m_1\op\m_2$ in
which $\m_1$ is a direct sum of finite dimensional irreducible
nontrivial  $\fg-$submodules and $\m_2$ is a trivial
$\fg-$submodule of $\LL.$ One can derive a specific vector space
$\fb$ from the $\fg-$module structure of $\m_1.$     This vector
space is equipped with  an algebraic structure which is induced by
the Lie algebraic structure of $\LL.$ Moreover the Lie algebra
$\LL$ can be reconstructed  from the algebra  $\fb$ in a
prescribed way, see \cite{ABG1} and \cite{ABG2}. 
This construction led to finding  a finite presentation for the universal central extension
of a Lie torus of a finite type other than $A$ and $C,$ see \cite{You}, \cite{AYY}.  This motivates us  to generalize this construction for Lie algebras graded by infinite root systems.

We give a complete description of the structure of root graded Lie algebras.
We fix an infinite irreducible locally finite  root system $R$ and show that a Lie algebra  $\LL$ graded by $R$ can be described in terms of a locally finite split simple Lie subalgebra   $\gg,$ some natural representations of $\gg$ and a certain  algebra  called the {\it coordinate algebra}. We also give the Lie bracket on $\LL$ in terms of the Lie bracket on $\gg,$ the action of the representations and the product on $\fb.$
More precisely, depending on type of $R,$ we consider a quadruple $\mathfrak{c}$
so called {\it coordinate quadruple}. We next correspond to
$\mathfrak{c},$ a specific algebra $\fb_{\mathfrak{c}}$ and a
specific Lie algebra $\{\fb_{\mathfrak{c}},\fb_{\mathfrak{c}}\}.$
Then for each subspace $\kk$ of the center of
$\{\fb_{\mathfrak{c}},\fb_{\mathfrak{c}}\}$ satisfying a certain
property called {\it the uniform property}, we define a Lie algebra
$\LL(\fb_{\mathfrak{c}},\kk)$ and show that it is a Lie algebra
graded by $R.$ Conversely, given a  Lie algebra $\LL$ graded
by $R,$ we prove that $\LL$ can be decomposed as $\m_1\op\m_2$
where $\m_1$ is a direct sum of  certain irreducible nontrivial
$\fg-$submodules for a locally finite spilt simple Lie subalgebra
$\fg$ of $\LL$  and $\m_2$ is a specific subalgebra of $\LL.$  We
derive a  quadruple $\mathfrak{c}$ from the $\fg-$module
structure of $\m_1$ and show that it is a coordinate quadruple. We
also prove that there is a subspace $\kk$ of
$\{\fb_{\mathfrak{c}},\fb_{\mathfrak{c}}\}$ satisfying the uniform
property such that $\m_2$ is isomorphic to the quotient algebra
$\{\fb_{\mathfrak{c}},\fb_{\mathfrak{c}}\}/\kk$ and moreover $\LL$
is isomorphic to $\LL(\fb_{\mathfrak{c}},\kk).$  If the root system $R$ is reduced, our method  suggests another approach to characterize 
Lie algebras graded by $R$ compared with what  is offered   by E.
Neher \cite{N}.

The author wishes to thank  the hospitality of Mathematics and
Statistics Department, University of Ottawa, where some parts of
this work were  carried out. The author also  would like to
express her sincere gratitude to Professor Saeid Azam and
Professor Erhard Neher for some fruitful discussions.

\section{Preliminaries}
Throughout   this work, $\bbbn$ denotes the set of  nonnegative
integers and $\bbbf$ is a field of characteristic zero. Unless
otherwise mentioned, all vector spaces are considered over
$\bbbf.$ We denote the dual space of a
vector space $V$ by $V^*.$ For a linear transformation $T$ on a vector space $V,$ if the trace of $T$ is defined, we denote it by  $tr(T).$ Also for a nonempty set $S,$ by $id_{_S}$ (or $id$ if there is no confusion), we mean the identity map on $S$ and by $|S|,$ we mean the cardinality of $S.$
Finally for an index set   $I,$ by a conventional notation, we take  $\bar
I:=\{\bar i\mid i\in I\}$ to be  a disjoint  copy of $I$ and for each subset $J$ of $I,$ by $\bar J,$ we mean the subset of $\bar I$ corresponding to $J.$

\subsection{Locally Finite Split Simple Lie Algebras}
In this subsection, we   recall the structure of  infinite dimensional
locally finite  split simple Lie algebras from \cite{NS} and state
some facts  which play   key roles in this work. Let us start with the following definition.
\begin{Definition}
{\em Let $\hh$ be a Lie algebra. We say an $\hh-$module $\m$ has a
{\it weight space decomposition with respect to $\hh,$} if
\begin{equation*}
\m=\op_{\a\in \hh^*}\m_\a \;\hbox{where}\; \m_\a:=\{x\in\m\mid
h\cdot x=\a(h)x;\;\;\forall\; h\in \hh\}
\end{equation*}
for all $\a\in\hh^*.$ The set $R:=\{\a\in\hh^*\mid \m_\a\neq
\{0\}\}$ is called the {\it set of weights} of $\m$ (with respect
to $\hh$).  For $\a\in R,$ $\m_\a$ is called a {\it weight space},
and  any element of $\m_\a$ is called a {\it weight vector} of
{\it weight} $\a.$  If a Lie algebra $\LL$ has a weight space
decomposition with respect to a nontrivial subalgebra $H$ of $\LL$
via the adjoint representation, $H$ is called a {\it split toral
subalgebra}. The set of weights of $\LL$ is called the {\it root
system} of $\LL$ with respect to $H,$ and  the corresponding
weight spaces are called {\it root spaces of $\LL.$} A Lie algebra
$\LL$ is called {\it split} if it contains a {\it splitting Cartan
subalgebra}, that is a  split toral subalgebra $H$ of $\LL$ with
$\LL_0=H.$ }
\end{Definition}

The root system of a locally finite split simple Lie algebra with
respect to a splitting Cartan subalgebra is a reduced irreducible locally finite root
system in the following sense (see \cite{Bo} and \cite{NS}):
\begin{Definition}\label{def-root}\cite{LN}
{\rm Let $\u$ be a nontrivial vector space and $R$ be a subset of
$\u,$ $R$ is said to be a {\it locally finite root system in $\u$}
of {\it rank} $dim(\u)$ if the following are satisfied:

(i) $R$ is locally finite, contains zero and spans $\u.$

(ii) For every $\a\in R^\times:=R\setminus\{0\},$ there exists
$\check\a\in\u^*$ such that $\check\a(\a)=2$ and $s_\a(\b)\in R$
for $\a,\b\in  R$ where $s_\a:\u\longrightarrow\u$ maps $u\in\u$
to $u-\check\a(u)\a.$ We set by convention $\check0$ to be zero.

(iii) $\check\a(\b)\in\bbbz,$ for $\a,\b\in R.$

Set  $R_{sdiv}:=(R\setminus\{\a\in R\mid 2\a\in R\})\cup \{0\}$
and call it the {\it semi-divisible  subsystem} of $R.$ The root
system $R$ is called {\it reduced} if $R=R_{sdiv}.$ }
\end{Definition}

Suppose that $R$  is a locally finite root   system. A nonempty
subset $S$ of $R$  is said to be  a {\it subsystem} of $R$ if $S$
contains zero and  $s_\a(\b)\in S$ for $\a,\b\in S\setminus\{0\}.$
A subsystem $S$ of $R$ is called {\it full} if $\hbox{span}_\bbbf
S\cap R=S.$ Following \cite[\S 2.6]{LN}, we say two nonzero roots
$\a,\b$ are {\it connected} if there exist finitely many roots
$\a_1=\a,\a_2,\ldots,\a_n=\b$ such that
$\check\a_{i+1}(\a_{i})\neq 0,$ $1\leq i\leq n-1.$ Connectedness
defines an equivalence relation on $R^\times$ and so $R^\times$ is the
disjoint union of its equivalence classes called {\it connected
components} of $R.$ A nonempty subset $X$ of $R$ is called {\it
irreducible,} if each two nonzero elements $x,y\in X$ are
connected and it is called {\it closed} if $(X+X)\cap R\sub X.$
 It is easy to see that
if $X$ is a connected component of a locally finite root system
$R,$ then $X\cup\{0\}$ is a closed subsystem of $R.$ For
the locally finite root system $R,$ take $\{R_\lam\mid
\lam\in\Gamma\}$ to be the class of all finite subsystems of $R,$
and say $\lam\preccurlyeq\mu$ $(\lam,\mu\in\Gamma)$ if $R_\lam$ is
a subsystem of $R_\mu,$ then $(\Gamma,\preccurlyeq)$ is a directed
set and $R$ is the direct union of $\{R_\lam\mid \lam\in\Gamma\}.$
Furthermore, if $R$ is irreducible, it is the direct union of its
irreducible finite subsystems.

Two locally finite root systems $( R,\u)$
 and $(S,\v)$ are  said  to be isomorphic if there is a linear transformation
$f:\u\longrightarrow \v$ such that $f(R)= S.$

Suppose that $I$ is a nonempty  index set  and $\u:=\op_{i\in
I}\bbbf\ep_i$ is the free $\bbbf-$module over   the
set $I.$ Define the  form $$\begin{array}{c}\fm:\u\times\u\longrightarrow\bbbf\\
(\ep_i,\ep_j)=\d_{i,j}, \hbox{ for } i,j\in I
\end{array}$$
and set
\begin{equation}\label{locally-finite}
\begin{array}{l}
\dot A_I:=\{\ep_i-\ep_j\mid i,j\in I\},\\
D_I:=\dot A_I\cup\{\pm(\ep_i+\ep_j)\mid i,j\in I,\;i\neq j\},\\
B_I:=D_I\cup\{\pm\ep_i\mid i\in I\},\\
C_I:=D_I\cup\{\pm2\ep_i\mid i\in I\},\\
BC_I:=B_I\cup C_I.
\end{array}
\end{equation}
One can see that these are irreducible locally finite root systems
in their $\bbbf-$span's which we refer to as  {\it type} $A,B,C,D$
and $BC$  respectively. Moreover every irreducible locally finite
root system of infinite rank is isomorphic to one of these root
systems (see \cite[\S4.14 
\S8]{LN}). Now we suppose $R$ is  an irreducible locally finite
root system as above and note that  $(\a,\a)\in\bbbn$ for all
$\a\in R.$ This allows us to  define
$$\begin{array}{l}
R_{sh}:=\{\a\in R^\times\mid (\a,\a)\leq(\b,\b);\;\;\hbox{for all $\b\in R$} \},\\
R_{ex}:=R\cap2 R_{sh},\\
R_{lg}:= R^\times\setminus( R_{sh}\cup R_{ex}).
\end{array}$$
The elements of $R_{sh}$ (resp. $R_{lg},R_{ex}$) are called {\it
short roots} (resp. {\it long roots, extra-long roots}) of $R$.

\medskip

 A locally finite split simple Lie algebra is said to be of type
$A,B,C$ or $D$ if its corresponding root system with respect to a
splitting Cartan subalgebra is of type $A,B,C$ or $D$
respectively. In what follows, we  recall from \cite{NS} the classification of infinite dimensional  locally finite split
simple Lie algebras.
 Suppose that $J$ is an index set and $\v=\v_J$ is a
vector space with a fixed basis $\{v_j\mid j\in J\}.$ One knows
that $\mathfrak{gl}(\v):=End(\v)$ together with
$$[\cdot,\cdot]:\mathfrak{gl}(\v)\times\mathfrak{gl}(\v)\longrightarrow
\mathfrak{gl}(\v);\; (X,Y)\mapsto XY-YX;\;\;
X,Y\in\mathfrak{gl}(\v)$$is a Lie algebra. Now for $j,k\in J,$
define
\begin{equation}\label{elementary2}e_{j,k}:\v\longrightarrow\v;\;\;
v_i\mapsto \d_{k,i}v_j,\;\;\; (i\in J),\end{equation} then
$\mathfrak{gl}(J):=\hbox{span}_\bbbf\{e_{j,k}\mid j,k\in J\}$ is a
Lie subalgebra of $\mathfrak{gl}(\v)$.

\begin{Lemma}[{Classical Lie algebras of type $A$}]\label{type-a-alg}
Suppose that $I$ is a non-empty index set of cardinality greater than $1$, $I_0$ is a fixed subset of
$I$ with $|I_0|>1$ and $\v$ is a vector space with a basis $\{v_i\mid i\in I\}.$
Take $\Lam$ to be an index set containing $0$ and  $\{I_\lam\mid \lam\in \Lam\}$ to be the class of all finite
subsets of $I$ containing $I_0.$ Set
$$\gg:=\mathfrak{sl}(I):=\{\phi\in\mathfrak{gl}(I)\mid
tr(\phi)=0\},$$ and  for $\lam\in\Lam,$   take
\begin{equation*}\label{simple-a-alg}
\begin{array}{l}\gg_{_{I_\lam}}:=\gg^\lam:=\gg\cap\hbox{span}\{e_{r,s}\mid r,s\in  I_\lam\}. \end{array}\end{equation*}
  Then
$\mathfrak{sl}(I)$ is a locally finite split simple Lie subalgebra
of $\mathfrak{gl}(I)$ with splitting Cartan subalgebra
$\hh:=\hbox{span}\{e_{i,i}-e_{j,j}\mid i,j\in I\}$ and
corresponding root system isomorphic to $\dot A_I.$ Moreover
for $i,j\in I$ with $i\neq j,$ we have $$\gg_{\ep_i-\ep_j}=\bbbf e_{i,j}.$$
Also for
each $\lam\in \Lam,$ $\gg^\lam$ is a finite dimensional split
simple Lie subalgebra of $\gg$ with  splitting Cartan subalgebra
$\hh^\lam:=\hh\cap\gg^\lam$ and  $\gg$ is the direct union of
$\{\gg^\lam\mid \lam \in \Lam\}.$

\end{Lemma}

In the following lemma, we see that  locally finite split simple
Lie algebras of type $B$ can be described in terms of derivations
of {\it Clifford Jordan algebras} which are  defined as   following:
\begin{Definition}[\cite{Yos}]\label{yoshii}
 {\rm Suppose that $A$ is a unital commutative associative algebra
over $\bbbf$ and $\w$ is a unitary $A$-module. Suppose that
$g:\w\times \w\longrightarrow A$ is a symmetric $A$-bilinear form and set
 $\jj=\jj(g,\w):=A\op\w.$ The vector space $\jj$ together with the following
multiplication
$$(a_1+w_1)(a_2+w_2)=a_1a_2+g(w_1,w_2)+a_1w_2+a_2w_1$$ for
$a_1,a_2\in A$ and $w_1,w_2\in\w$ is a Jordan algebra called a
{\it Clifford Jordan algebra.} For $a,b\in\jj,$ define
$D_{a,b}:=-[{\bf L}_a,{\bf L}_b]:={\bf L}_b{\bf L}_a-{\bf L}_a{\bf L}_b$ where ${\bf L}_a,{\bf L}_b$ are
left multiplications by $a$ and  $b$ respectively. For a subspace
$V$ of $\jj,$ set  $D_{V,V}$ to be the subspace of endomorphisms
of $\jj$ spanned by $D_{a,b}$ for $a,b\in V.$ One can see that for
$w_1,w_2\in\w,$ $D_{w_1,w_2}$ can be identified with
$D_{w_1,w_2}|_{_\w}.$ This  allows us to consider $D_{\w,\w}$ as a
subalgebra of $\mathfrak{gl}(\w).$}
\end{Definition}

\begin{Lemma}[{Classical Lie algebras of type $B$}]

\label{type-b-alg} Suppose that $I$ is a non-empty index set. Take
$J:=\{0\}\uplus I\uplus\bar{I}$ and consider the vector space
$\v:=\v_J$ as before. Define the bilinear form $\fm$ on $\v$ by
\begin{equation}\label{form-b-alg}\begin{array}{c}(v_j,v_{\bar k})=(v_{\bar k},v_j)=2\d_{j,k},\; (v_0,v_0)=2,\\
(v_j,v_k)=(v_j,v_0)=(v_0,v_j)=(v_0,v_{\bar j})=(v_{\bar
j},v_0)=(v_{\bar j},v_{\bar k})=0;\;j,k\in
I,\end{array}\end{equation} and  set
$$\gg:=\mathfrak{o}_B(I):=\{\phi\in\mathfrak{gl}(J)\mid
(\phi(v),w)=-(v,\phi(w)),\; \hbox{for all $v,w\in\v$}\}.$$
Then we have the following:

(i) $\gg$ is a  locally finite  split simple Lie subalgebra of
$\mathfrak{gl}(J)$ with splitting Cartan subalgebra
$\hh:=\hbox{span}_\bbbf\{h_i:=e_{i,i}-e_{\bar i,\bar i}\mid i\in
I\}$  and corresponding root system isomorphic to $B_I.$ Moreover for $i,j\in J$ with $i\neq j,$ we have
$$\begin{array}{c}\gg_{\ep_i-\ep_j}=\bbbf (e_{i,j}-e_{\bar j,\bar i}),\;\gg_{\ep_i+\ep_j}=\bbbf (e_{i,\bar j}-e_{ j,\bar i}),\;\gg_{-\ep_i-\ep_j}=\bbbf (e_{\bar i,j}-e_{\bar j,i})\\
\gg_{\ep_i}=\bbbf (e_{i,0}-e_{0,\bar i}),\;\gg_{-\ep_i}=\bbbf (e_{\bar i,0}-e_{0,i}).\end{array}$$

(ii) For  the Clifford Jordan algebra $\jj(\fm,\v),$ we have  $\gg=D_{\v,\v}.$

(iii) For a fixed subset $I_0$ of $I,$ take $\Lam$ to be an index set containing $0$ such that  $\{I_\lam\mid
\lam\in\Lam\}$ is the class of all finite subsets of $I$
containing $I_0.$ For each $\lam\in \Lam,$ set
\begin{equation}\label{simple-b-alg}\begin{array}{l}\gg_{_{I_\lam}}:=\gg^\lam:=\gg\cap\hbox{span}\{e_{r,s}\mid r,s\in \{0\}\cup I_\lam\cup\bar
I_\lam\}. \end{array}\end{equation} Then $\gg^\lam$
($\lam\in\Lam$) is a finite dimensional split simple Lie
subalgebra of $\gg$ of type $B,$ with  splitting Cartan subalgebra
$\hh^\lam:=\hh\cap\gg^\lam$ and  $\gg$ is the  direct union of
$\{\gg^\lam\mid \lam \in \Lam\}.$

\end{Lemma}
\begin{Lemma}
[Classical Lie algebras of type $D$] \label{type-d-alg} Suppose
that $I$ is a non-empty index set and $I_0$ is a fixed subset of
$I.$ Set $J:=I\uplus\bar{I}$  and take $\{I_\lam\mid
\lam\in\Lam\},$ where $\Lam$ is an index set containing $0,$ to be the class of all finite subsets of $I$
containing $I_0.$ Define the bilinear form $\fm$ on $\v=\v_J$ by
\begin{equation}\label{form-d-alg}(v_j,v_{\bar k})=(v_{\bar k},v_j)=2\d_{j,k},\;
(v_j,v_k)=(v_{\bar j},v_{\bar k})=0;\;(j,k\in I),\end{equation}
and set
$$\begin{array}{l}\gg:=\mathfrak{o}_D(I):=\{\phi\in\mathfrak{gl}(J)\mid
(\phi(v),w)=-(v,\phi(w)),\; \hbox{for all $v,w\in\v$}\},\\
\hh:=\hbox{span}_\bbbf\{h_i:=e_{i,i}-e_{\bar i,\bar i}\mid i\in
I\}.\end{array}
$$Also for $\lam\in\Lam,$ take $$\gg_{_{I_\lam}}:=\gg^\lam:=\gg\cap\hbox{span}\{e_{r,s}\mid r,s\in I_\lam\cup\bar
I_\lam\}.$$  Then $\gg$ is a  locally finite  split
simple Lie subalgebra of $\mathfrak{gl}(J)$ with splitting Cartan
subalgebra $\hh$ and corresponding root system isomorphic to
$D_I.$ Moreover for $i,j\in J$ with $i\neq j,$ we have
$$\gg_{\ep_i-\ep_j}=\bbbf (e_{i,j}-e_{\bar j,\bar i}),\;\gg_{\ep_i+\ep_j}=\bbbf (e_{i,\bar j}-e_{ j,\bar i}),\;\gg_{-\ep_i-\ep_j}=\bbbf (e_{\bar i,j}-e_{\bar j,i}).$$ Also for each $\lam\in \Lam,$ $\gg^\lam$ is a finite
dimensional split simple Lie subalgebra of $\gg,$ of type $D,$
with  splitting Cartan subalgebra $\hh^\lam:=\hh\cap\gg^\lam,$ and
$\gg$ is the direct union of $\{\gg^\lam\mid \lam \in \Lam\}.$

\end{Lemma}

\begin{Lemma}
[Classical Lie algebras of type $C$]\label{type-c-alg}  Suppose
that $I$ is a non-empty index set and $J:= I\uplus\bar{I}.$
Consider the bilinear form $\fm$ on $\v=\v_J$ defined
by\begin{equation}\label{form-c}(v_j,v_{\bar k})=-(v_{\bar
k},v_j)=2\d_{j,k},\; (v_j,v_k)=0,\; (v_{\bar j},v_{\bar
k})=0,\;\;\;  (j,k\in I),\end{equation}and set
$$\gg:=\mathfrak{sp}(I):=\{\phi\in\mathfrak{gl}(J)\mid
(\phi(v),w)=-(v,\phi(w)),\; \hbox{for all $v,w\in\v$}\}.$$ Also for a fixed subset $I_0$ of $I,$ take $\{I_\lam\mid
\lam\in\Lam\}$ to be the class of all finite subsets of $I$
containing $I_0,$ in which $\Lam$ is an index set containing $0,$ and  for each $\lam\in \Lam,$ set
\begin{equation}\label{simple-c-alg}\gg_{_{I_\lam}}:=\gg^\lam:=\gg\cap\hbox{span}\{e_{r,s}\mid r,s\in I_\lam\cup\bar
I_\lam\}.
\end{equation}
 Then $\gg$
 is a  locally finite  split simple Lie
subalgebra of $\mathfrak{gl}(J)$ with splitting Cartan subalgebra
$\hh:=\hbox{span}_\bbbf\{h_i:=e_{i,i}-e_{\bar i,\bar i}\mid i\in
I\}.$ Moreover for $i,j\in I$ with $i\neq j,$ we have
$$\begin{array}{c}\gg_{\ep_i-\ep_j}=\bbbf (e_{i,j}-e_{\bar j,\bar i}),\;\gg_{\ep_i+\ep_j}=\bbbf (e_{i,\bar j}+e_{ j,\bar i}),\;\gg_{-\ep_i-\ep_j}=\bbbf (e_{\bar i,j}+e_{\bar j,i})\\
\gg_{2\ep_i}=\bbbf e_{i,\bar i},\;\gg_{-2\ep_i}=\bbbf e_{\bar i,i}.\end{array}$$ Also for $\lam\in\Lam,$ $\gg^\lam$ is  a finite dimensional
split simple Lie subalgebra  of type $C,$ with  splitting Cartan
subalgebra $\hh^\lam:=\hh\cap\gg^\lam,$ and
 $\gg$ is the
direct union of $\{\gg^\lam\mid \lam\in\Lam\}.$
\end{Lemma}
\begin{Proposition}\label{locally-alg}\cite[Theorem VI.7]{NS}
Suppose that $I$ is an infinite index set, then $\mathfrak{o}_B(I)$ is isomorphic to $\mathfrak{o}_D(I).$ Moreover if $\gg$ is an infinite dimensional  locally finite split simple Lie algebra, then $\gg$ is isomorphic to exactly one of the Lie algebras $\mathfrak{sl}(I),$ $\mathfrak{o}_B(I)$ or $\mathfrak{sp}(I).$
\end{Proposition}

\begin{Lemma}\label{final3}
Suppose that $R$ is an irreducible locally finite root system and
$S$  is an irreducible closed subsystem of $R.$ Suppose   $\gg$ is
a locally finite split simple Lie algebra with a splitting Cartan
subalgebra $\hh$ and the root system $R_{sdiv}.$ Set
$\fg:=\sum_{\a\in S_{sdiv}^\times}\gg_\a\op\sum_{\a\in
S_{sdiv}^\times}[\gg_\a,\gg_{-\a}]$ and $\fh:=\hh\cap\fg,$  then
the restriction of
$$\pi:\hh^*\longrightarrow \fh^*;\;\;f\mapsto f|_{\fh},\;\;
f\in \hh^*$$ to  $S$ is injective. Identify $\a\in S$ with
$\pi(\a)$ via $\pi,$  then   $\fg$ is a locally finite split simple Lie
subalgebra of $\gg$ with  splitting Cartan subalgebra $\fh$ and corresponding
root system $S_{sdiv}.$
\end{Lemma}

\pf  We first  claim that
\begin{equation}\label{crutial}\parbox{3in}{\it\begin{center} if  $\a,\b\in S$ and   $\a-\b\not \in R,$
then  there is $h\in \fh$ such that $\a(h)>0$ and
$\b(h)\leq0.$\end{center}}\end{equation} To prove this, we note
that  since $\a-\b\not\in R,$ we have  $\a\neq 0$ and $\b\neq0 .$
Moreover,  it follows  from the theory of locally finite root
systems  that  $\b-2\a\not\in R$ and  $\a-2\b\not\in R,$ also if
$2\a\in R$ or $2\b\in R,$ then $2\a-2\b\not \in R.$ Therefore
setting
$$\a':=\left\{\begin{array}{ll}\a&\hbox{ if $2\a\not\in R$}
\\
2\a&\hbox{ if $2\a\in R,$}
\end{array}\right.\andd \b':=\left\{\begin{array}{ll}\b&\hbox{ if $2\b\not\in R$}
\\
2\b&\hbox{ if $2\b\in R,$}\end{array}\right. $$ we have
$\a',\b'\in S_{sdiv}^\times$ and $\a'-\b'\not\in R.$ Next we fix
$e\in\gg_{\a'}$ and $f\in\gg_{-\a'}$ such that $(e,h:=[e,f],f)$ is
an $\mathfrak{sl}_2$-triple. Since $\a'-\b'\not\in R_{sdiv},$ one
knows from $\mathfrak{sl}_2-$module theory  that $\b'(h)\leq 0$
while $\a'(h)=2>0.$ Therefore $h\in[\gg_{\a'},\gg_{-\a'}]\sub\fh,$
$\a(h)>0$ and $\b(h)\leq 0.$ 
This completes the proof of the
claim. Now suppose $\a,\b\in S$ with $\pi(\a)=\pi(\b).$ We must
show $\a=\b.$ We prove this through the  following three cases:

\underline{\bf Case 1. $\a,\b\in S_{sdiv}:$} If  $\gamma:=\a-\b\in
R^\times,$ then since $S$ is a closed subsystem of $R$  and
$S_{sdiv}$ is a closed subsystem of $S,$  we get $\gamma\in
S_{sdiv}^\times.$ Thus there is
$t\in[\gg_\gamma,\gg_{-\gamma}]\sub\fh$ with $\gamma(t)=2,$ so
$(\a-\b)(t)=2$ which contradicts the fact that
$\a\mid_{\fh}=\b\mid_{\fh}.$ Therefore  $\a-\b\not\in R^\times. $
Now if $\a-\b\neq0,$ then $\a-\b\not\in R$ and so  using
(\ref{crutial}), one  finds $h\in\fh$ with $\a(h)>0$ and
$\b(h)\leq0.$ This  is again a contradiction. Therefore $\a=\b.$

 \underline{\bf Case 2. $\a,\b\not\in S_{sdiv}:$} In this case $2\a,2\b\in S_{sdiv}$ and so by Case 1, $2\a=2\b$ which in turn implies that $\a=\b.$

\smallskip

\underline{\bf Case 3. $\a\in S_{sdiv}, \b\not\in S_{sdiv}:$}  If
$\a-\b\not\in R,$ then by (\ref{crutial}), there is $h\in \fh$
such that $\a(h)>0$ and $\b(h)\leq 0$ which contradicts the fact
that $\pi(\a)=\pi(\b).$ Also if $\a=2\b,$ then since $\a\in
S_{sdiv}^\times,$ there is $h\in[\gg_\a,\gg_{-\a}]\sub\fh$ with
$\a(h)=2.$ Thus $\a(h)\neq \b(h)$ which is again  a
contradiction. Therefore $\a-\b\in R$ and $\a\neq 2\b.$ Now  if
$\a-\b\neq0,$ we get that  $\a-\b\in R_{sh},$  $\a\in R_{lg},$
$\b\in R_{sh},$ $\gamma:=\a-2\b\in  R_{lg},$
$\a+\gamma,\a-\gamma\in R_{ex}$ and $\a+2\gamma,\a-2\gamma\not\in
R.$ Now since $\gamma\in S_{sdiv}^\times,$ there is
$h\in[\gg_\gamma,\gg_{-\gamma}]\sub\fh$ with $\gamma(h)=2.$ Also
since $\a+2\gamma,\a-2\gamma\not\in R,$ one concludes form
$\mathfrak{sl}_2-$module theory that $\a(h)=0.$
So we have $\b(h)=\a(h)=0.$ But this
gives that $2=\gamma(h)=(\a-2\b)(h)=0,$ a contradiction. Thus we
have  $\a=\b.$ This completes the proof of the first assertion.

For the last assertion, we note that  $S$ is a closed subsystem of $R$ and $S_{sdiv}$
is  a closed subsystem of $S,$ so it is easily seen that $\fg$ is
a subalgebra of $\gg.$ Now this together with the fact that
$\pi|_S$ is injective completes the proof.\qed

\bigskip

\begin{Definition}\label{direct-limit}
{\rm Take  $\gg,$ $\Lam,$ $\gg^\lam,$ and $\hh^\lam$ $(\lam\in\Lam)$ to
be as in one of Lemmas \ref{type-a-alg}, \ref{type-b-alg},
\ref{type-d-alg} and \ref{type-c-alg}. For $\lam,\mu\in\Lam,$ we
say $\lam\preccurlyeq\mu$ if $\gg^\lam$ is a subalgebra of
$\gg^\mu.$ Let $\chi$ be a representation of $\gg$ in a vector
space $\m.$ We say $\m$ is a {\it direct limit $\gg$-module with
directed system} $\{\m^\lam\mid \lam\in\Lam\}$ if
\begin{itemize}
\item  for $\lam,\mu\in\Lam$ with $\lam\preccurlyeq\mu,$
$\m^\lam\sub\m^\mu\sub\m$ and as a vector space, $\m$ is the direct
union of $\{\m^\lam\mid\lam\in\Lam\},$ \item for $\lam\in\Lam,$
$\m^\lam$ is finite dimensional and for all $x\in\gg^\lam,$
$\m^\lam$ is invariant under $\chi(x),$ \item for $\lam\in\Lam,$
$\chi\mid_{\gg^\lam}$ defines a  nontrivial finite dimensional
irreducible  $\gg^\lam-$module in $\m^\lam$ having a weight space decomposition with respect to $\hh^\lam$ whose  weight spaces
corresponding to nonzero weights are one dimensional, \item for
$\lam,\mu\in\Lam$ with $\lam\preccurlyeq\mu,$  the set of weights
of $\gg^\lam-$module  $\m^\lam$ is contained in the set of weights of $\gg^\mu-$module  $\m^\mu$
restricted to $\hh^\lam$ and
$(\m^\lam)_{p|_{_{\hh^\lam}}}=(\m^\mu)_p$ for each  weight $p$ of $\m^\mu$ for which $p|_{_{\hh^\lam}}$ is a  nonzero weight
 of $\m^\lam.$
\end{itemize}}
\end{Definition}
Using standard techniques, one can verify the following propositions:
\begin{Proposition}\label{dir-lim-mod} Consider $\gg,$ $\Lam,$  $\gg^\lam$ and  $\hh^\lam$  $(\lam\in\Lam)$ as in Definition \ref{direct-limit}. Suppose that
$\m$ is a  direct limit $\gg-$module with directed system
$\{\m^\lam\mid\lam\in\Lam\}.$ Then we have the followings:

(i) $\m$ is an irreducible $\gg-$module.

(ii) If  $\w$ is another  direct limit $\gg-$module
with directed system $\{\w^\lam\mid \lam\in\Lam\},$ and  for each
$\lam\in\Lam,$ $\gg^\lam-$modules  $\m^\lam$ and $\w^\lam$ are
isomorphic, then as two $\gg-$modules, $\m$ and  $\w$ are
isomorphic.
\end{Proposition}
\medskip

\begin{Proposition}
\label{rep-local} Suppose that $I$ is a nonempty index set.

(a)  Take $\gg:=\mathfrak{o}_B(I)$ and use the same notations as in Lemma \ref{type-b-alg}. Define $$\pi:\gg\longrightarrow
End(\v);\;\pi(\phi)(v)=\phi(v);\;\;\phi\in\gg,\; v\in \v,$$ then

(i) $\pi$ is an irreducible  representation of $\gg$ in $\v$ equipped with a weight space decomposition with respect to $\hh$ whose set
of weights is $\{0,\pm\ep_i\mid i\in I\}$ with $\v_0=\bbbf v_0,$
$\v_{\ep_i}=\bbbf v_i$ and $\v_{-\ep_i}=\bbbf v_{\bar i}$ for
$i\in I,$

(ii) for  each $\lam\in \Lam,$ set
\begin{equation}\label{simple-b-mod}\v_{_{I_\lam}}:=\v^\lam:=\hbox{span}_\bbbf\{v_r\mid r\in \{0\}\cup
I_\lam\cup\bar I_\lam\},\end{equation} then $\v$ is a direct limit $\gg-$module with directed system
$\{\v^\lam\mid \lam \in \Lam\}.$

(b) Use the same notations as in  Lemma \ref{type-c-alg} and take  $\gg:=\mathfrak{sp}(I).$

(i) Define
$$\pi_1:\gg\longrightarrow End(\v);\;\pi(\phi)(v):=\phi(v);\;\;\phi\in\gg,\; v\in
\v.$$ Then $\pi_1$ is an irreducible  representation of $\gg$ in
$\v$ equipped with a weight space decomposition with respect to $\hh$ whose set of weights is $\{\pm\ep_i\mid i\in I\}$ with
$\v_{\ep_i}=\bbbf v_i$ and $\v_{-\ep_i}=\bbbf v_{\bar i}$ for
$i\in I.$ Also for $J:=I\cup\bar{I}$ and
\begin{equation}\label{module-s-c}\ss:=\{\phi\in \mathfrak{gl}(J)\mid
tr(\phi)=0,(\phi(v),w)=(v,\phi(w)),\; \hbox{for all
$v,w\in\v$}\},\end{equation} $$\pi_2:\gg\longrightarrow
End(\ss);\; \pi_2(X)(Y):=[X,Y];\;\; X\in \gg,\; Y\in \ss$$ is an
irreducible representation of $\gg$ in $\ss$  equipped with a weight space decomposition with respect to $\hh$ whose set of weights
is $\{0,\pm(\ep_i\pm\ep_j)\mid i,j\in I,\; i\neq j\}$ with
$\ss_0=\hbox{span}_\bbbf\{e_{r,r}+e_{\bar r,\bar
r}-\frac{1}{|I_\lam|}\sum_{i\in I_\lam} (e_{i,i}+e_{\bar i,\bar i})\mid
\lam\in\Lam,r\in I_\lam\},$ $\ss_{\ep_i+\ep_j}=\bbbf(e_{i,\bar
j}-e_{j,\bar i}),$ $\ss_{-\ep_i-\ep_j}=\bbbf(e_{\bar i,j}-e_{\bar
j, i})$ and $\ss_{\ep_i-\ep_j}=\bbbf(e_{i,j}+e_{\bar j,\bar i})$
($i,j\in I, i\neq j$).

(ii) For  $\lam\in \Lam,$ set
\begin{equation}\label{simple-c}\begin{array}{l}
\v_{_{I_\lam}}:=\v^\lam:=\hbox{span}_\bbbf\{v_r\mid r\in I_\lam\cup\bar
I_\lam\},\\\ss_{_{I_\lam}}:=\ss^\lam:=\ss\cap\hbox{span}_\bbbf\{e_{r,s}\mid r,s\in
I_\lam\cup\bar I_\lam\}.\end{array}\end{equation}  Then $\v$ is a  direct limit $\gg-$module with directed system
$\{\v^\lam\mid \lam \in \Lam\}$ and $\ss$ is a direct limit $\gg-$module
with directed system   $\{\ss^\lam\mid \lam\in \Lam\}.$
\end{Proposition}

\subsection{Finite Dimensional Case} In this subsection, we state a proposition on representation theory of finite dimensional  split simple Lie algebras. This proposition  is an essential tool for the proof of our results in  the next section. We start with an elementary but important fact  about finite dimensional representations of a finite dimensional split semisimple Lie algebra.
\begin{Lemma}\label{elementary}
Suppose that $\gg$ is a finite dimensional split semisimple Lie
algebra with a splitting Cartan subalgebra $\hh$ and the root system
$R.$ Let  $\v$ be   a finite dimensional $\gg$-module equipped with a weight space decomposition with respect to $\hh.$ Take $\Pi$ to
be  the set of weights of $\v.$ If $\a\in
R^\times$ and $\lam\in \Pi$ are such that $\a+\lam\in\Pi,$ then
$\gg_\a\cdot\v_\lam\neq\{0\}.$ In particular if $\v_{\lam+\a}$ is
one dimensional, then $\gg_\a\cdot\v_\lam=\v_{\lam+\a}.$
\end{Lemma}
\pf  Take $e\in\gg_\a$ and  $f\in\gg_{-\a}$ to be such that
$(e,h:=[e,f],f)$ is an $\mathfrak{sl}_2$-triple  and define
$\mathfrak{s}:=\hbox{span}_\bbbf\{e,h,f\}.$ Set
$\w:=\sum_{k=-\infty}^{\infty}\v_{\lam+k\a},$ then $\w$ is a
finite dimensional  $\mathfrak{s}$-module and so by Weyl theorem,
it is decomposed into finite dimensional irreducible
$\mathfrak{s}$-modules, say $\w=\op_{i=1}^n\w_i$ where $n$ is a
positive integer and  $\w_i,$ $1\leq i\leq n,$ is a finite
dimensional irreducible $\mathfrak{s}$-module. We note that the
set of weights of $\w$ with respect to $\bbbf h$ is
$\Pi'=\{\lam(h)+2k\mid k\in\bbbz\andd \lam+k\a\in\Pi\}$ and that
for $k\in\bbbz$ with $ \lam+k\a\in\Pi,$
$\w_{\lam(h)+2k}=\v_{\lam+k\a}.$ Now as $\lam,\lam+\a\in\Pi,$ we
have $\lam(h),\lam(h)+2\in\Pi'$ and so by $\mathfrak{sl}_2$-module
theory, there is $1\leq i\leq n$ such that $\lam(h),\lam(h)+2$ are
weights for $\w_i.$ Now again using $\mathfrak{sl}_2$-module
theory, we get that $$0\neq e\cdot (\w_i)_{\lam(h)}\sub e\cdot
\w_{\lam(h)}= e\cdot \v_\lam\sub\gg_\a\cdot\v_\lam$$ showing that $\gg_\a\cdot\v_\lam\neq\{0\}.$ The last
statement is  derived  simply from the first assertion.\qed

\begin{Lemma}\label{gen-fact}Suppose that $\{e_i,f_i,h_i\mid 1\leq i\leq n\}$ is a set of
Chevalley generators for a finite dimensional split simple Lie
algebra  $\gg$ and $\v$ is a $\gg-$module equipped with  a weight
space decomposition with respect to the Cartan subalgebra
$\hh:=\hbox{span}\{h_i\mid 1\leq i\leq n\}.$
 Let $v$ be a weight
vector, $m$ be a positive integer and  $1\leq i, j_1,\ldots,j_m\leq
n.$ Let the set $\{k\in\{1,\ldots,m\}\mid j_k=i\}$ be a
nonempty set and $k_1<\cdots<k_p$ be such that   $\{k\in\{1,\ldots,m\}\mid j_k=i\}=\{k_1,\ldots,k_p\}.$  Then if
$f_i\cdot v=0,$ we have 
$$f_i\cdot e_{j_m}\cdot\cdots\cdot e_{j_1}\cdot
v\in\sum_{t=1}^p\bbbf e_{j_m}\cdot\cdots e_{j_{k_p}}\cdot\cdots
\cdot \widehat{e_{j_{k_t}}}\cdot\cdots\cdot e_{j_{k_1}}\cdot\cdots
\cdot e_{j_1}\cdot v,$$ in which  $"\widehat{\;\;}\;"$ means
omission.
\end{Lemma}

\pf Using Induction on $p,$ we are done.

\qed

\begin{Proposition}
\label{gen1} Suppose $R_1$ is an irreducible finite root system
and $R_2$ is an irreducible full subsystem of $R_1$ of rank greater that 1.  Let
$\gg_1$ be  a finite dimensional split simple Lie algebra with a
splitting Cartan subalgebra $\hh_1$ and corresponding  root system
$(R_1)_{sdiv}.$ Set $\gg_2:=\sum_{\a\in
(R_2)_{sdiv}^\times}(\gg_1)_\a\op\sum_{\a\in
(R_2)_{sdiv}^\times}[(\gg_1)_\a,(\gg_1)_{-\a}]$ and
$\hh_2:=\hh_1\cap\gg_2.$ For $i=1,2,$ assume  $\v_i$ is a
$\gg_i$-module equipped with a weight space decomposition with
respect to $\hh_i$ and take $\Lam_i$ ($i=1,2$) to be the set of
wights of $\v_i$ with respect to $\hh_i.$ Suppose that

(i)  $R_1$ and $R_2$ are of the same type $X\neq
G_2,F_4,E_{6,7,8}$,

(ii) $\Lam_1\sub R_1$ and $\Lam_2\sub\{\a|_{_{\hh_2}}\mid \a\in
R_2\},$



(iii) $\v_2\sub\v_1$ with $(\v_2)_{\a|_{_{\hh_2}}}\sub(\v_1)_\a,$
for $\a\in \{\b\in R_2\mid\b|_{\hh_2}\in\Lam_2\}\setminus\{0\}.$

Let  $\w$ be a nontrivial finite dimensional irreducible
$\gg_2$-submodule of $\v_2$ and  take $\u$ to be the
$\gg_1$-submodule of $\v_1$ generated by $\w,$ then $\u$ is a
finite dimensional irreducible $\gg_1-$module equipped with a weight space decomposition with respect to $\hh_1$ whose set of nonzero
weights is $(R_1)_{sh},$  (resp. $(R_1)^\times_{sdiv},$ or
$((R_1)_{sdiv})_{sh}$) if the set of nonzero  weights of $\w$ is the set of elements of
$(R_2)_{sh}$ (resp. $(R_2)^\times_{sdiv},$ or $((R_2)_{sdiv})_{sh}$) restricted to $\hh_2.$
\end{Proposition}
\pf  Take  $n:=n_1$ and $\ell:=n_2$ to be the rank of $R_1$ and $R_2$
respectively. Using Lemma \ref{final3}, we identify $\b\in R_2$ with $\b|_{\hh_2}.$ Also without loss of generality, we assume  $R_1,R_2$ and
bases $\D _1,\D_2$ for $(R_1)_{sdiv}$ and $(R_2)_{sdiv}$
respectively are as in the following tables:$${\small
\begin{tabular}{|c|c|}
\hline
\hbox{ Type }& \hbox{   $R_k(k=1,2)$ }\\
\hline
$ A$& $\{\pm(\ve_i-\ve_j)\mid  1\leq i<j\leq n_k+1\}\cup\{0\},$ $n_k\geq2$\\
\hline
$B$&$\{\pm\ve_i,\pm(\ve_i\pm \ve_j)\mid  1\leq i<j\leq n_k\}\cup\{0\},$ $n_k\geq2$\\
\hline
$C$&$\{\pm2\ve_i,\pm(\ve_i\pm \ve_j)\mid  1\leq i<j\leq n_k\}\cup\{0\},$ $n_k\geq3$\\
\hline
$D$ &$\{\pm(\ve_i\pm \ve_j)\mid  1\leq i<j\leq n_k\}\cup\{0\},$ $n_k\geq4$\\
\hline
$BC$&$\{\pm\ve_i,\pm(\ve_i\pm \ve_j)\mid 1\leq i,j\leq  n_k\},$ $n_k\geq2$\\
\hline
\end{tabular}
}$$

$$  {\small
\begin{tabular}{|c|c|}
\hline
\hbox{ Type }& \hbox{ $\D_k(k=1,2)$ }\\
\hline
$ A$& $\{\a_i:=\ep_{i+1}-\ep_{i}\mid 1\leq i\leq n_k\}$\\
\hline
$B$&$\{\a_1:=\ep_1,\a_i:=\ep_i-\ep_{i-1}\mid 2\leq i\leq n_k\}$\\
\hline
$C$&$\{\a_1:=2\ep_1,\a_i:=\ep_i-\ep_{i-1}\mid 2\leq i\leq n_k\}$\\
\hline
$D$ &$\{\a_1:=\ep_1+\ep_2,\a_i:=\ep_i-\ep_{i-1}\mid 2\leq i\leq n_k\}$\\
\hline
$BC$&$\{\a_1:=2\ep_1,\a_i:=\ep_i-\ep_{i-1}\mid 2\leq i\leq n_k\}$\\
\hline
\end{tabular}
} $$

Suppose that $\w$ is  of highest pair
$(v,\a)$  with respect to $\D_2.$  Since   the set of
weights of $\w$ is permuted by the Weyl group of $R_2,$ one gets that
$\a=\a_{*}^2$ for $*\in\{sh,lg,ex\}$ where for $i=1,2,$ $\a^i_{sh}$ (resp. $\a^i_{lg},$ or $\a^i_{ex}$)  denotes the highest short (resp. long, or extra long) root
of $R_i$ with respect to $\D_i.$  Next suppose that  $\{e_i,f_i,h_i\mid
1\leq i\leq n\}$ is a set of Chevalley generators for $\gg_1$ with respect to $\D_1,$
then $\{e_i,f_i,h_i\mid 1\leq i\leq \ell\}$ is  a set of Chevalley
generators for $\gg_2$ with respect to $\D_2.$ Now  as  $\u$ is a $\gg_1-$submodule of $\v_1$ generated by $v,$ we have
\begin{equation}\label{tavakol1}\u=\sum_{t,s\in\bbbn}\bbbf (f_{i_t}\cdot\cdots\cdot
f_{i_1}\cdot e_{j_s}\cdot \cdots\cdot e_{j_1}\cdot v)\end{equation}
where $i_1,\ldots,i_t,j_1,\ldots,j_s\in\{1,\ldots,n\}.$ This implies that  $\u$ is finite dimensional as  $\Lam_1$ is a
subset of the finite root system $ R_1.$  So there are a positive integer $p$ and irreducible finite dimensional $\gg_1-$submodules $\u_j$ ($1\leq
j\leq p$) of $\u$ such that  $\u=\op_{j=1}^p\u_j.$
But we  know that $v$
generates the $\gg_1-$submodule  $ \u,$  and that  $v\in \u\cap(\v_2)_\a\sub\u\cap(\v_1)_\a=\u_\a=
\op_{j=1}^p(\u_j)_\a,$ so

\begin{equation}\label{gen3}
\parbox{2.6in}{for any  $1\leq j\leq p,$ there is a nonzero element $u_j\in(\u_j)_\a$ such that
$v=\sum_{j=1}^pu_j$.}
\end{equation}
This in particular  implies that each $\u_j$ $(1\leq j\leq p)$ is a  nontrivial irreducible  $\gg_1-$module.
But we know that  for $1\leq
j\leq p,$ the set of weights of $\u_j$ is a subset of $R_1,$ and that it is permuted by the Weyl group of $R_1,$ so the
highest weight of  $\u_j$ is $\a_{*}^1$ for $*=sh,lg,ex.$ Therefore using the   finite dimensional theory, one knows  that
\begin{equation}\label{one-dim}\parbox{4.5in}{  \begin{center}{\it the weight spaces of $\u_j$ ($1\leq j\leq p$ ) with respect to $\hh_1$ corresponding to nonzero weights are one dimensional.}\end{center}}
\end{equation}

Now we are ready to proceed  with the proof   in the following
three steps:
\begin{itemize}
\item{Step 1.} $ \u_t$ ($1\leq t\leq p$) is  a finite dimensional
irreducible $\gg_1-$module of highest weight $\a_*^1$ if
$\a=\a_*^2$ for $*=sh,lg,ex.$ \item{Step 2.} $\dim( \u_\a)=1,$
\item{Step 3.} $p=1.$
\end{itemize}

\noindent\textbf{Step 1:} We use a case-by-case argument to prove the desired point. We  note that there is nothing to show   if $R_1$ is of type
$A$ or $D,$ and    continue as following:
\medskip

\textbf{\emph{\underline{Type $B:$}}} One can see that in this
case $\a_{sh}^1=\ep_n,$ $\a_{sh}^2=\ep_\ell,$
$\a_{lg}^1=\ep_{n}+\ep_{n-1}$ and
$\a_{lg}^2=\ep_\ell+\ep_{\ell-1}.$ We  first assume $\a=\a^2_{sh}$
and    show that the highest weight of $\u_t$ is the highest short
root of $R_1$. For this, it is enough to show that no long root is
a weight for $\u_t.$ Suppose to the contrary that the set of
weights of $\u_t$ contains a long root or equivalently contains all long roots. Setting $\b:=\ep_{\ell-1},$ we get
$\a+\b$ is a long root  of $R_1$ and so $\a+\b$ is  a weight for
$\u_t.$ Now fix $x\in(\gg_2)_{\b}=(\gg_1)_\b$ and note that  $\a+\b$ is a weight for
$ \u_t;$ applying  Lemma \ref{elementary} together with (\ref{gen3})  and
(\ref{one-dim}), we have   $x\cdot u_t\neq 0.$ This gives that   $0\neq \sum_{j=1}^px\cdot u_j=x\cdot
v\in
(\gg_2)_\b\cdot\w_\a\sub \w_{\a+\b}$  which is a contradiction as $\a+\b$ is a long root and cannot be a weight for $\w.$ Therefore
$\u_t$ has no long root as a weight  and so we are done in the
case $\a=\a_{sh}^2$. Next suppose $\a=\a^2_{lg}$ and note that
by (\ref{gen3}), $u_t$ is a weight
vector of $\u_t$ of weight $\a.$ Since $\a$ is a long root, the set of  weights of  $\u_t$ contains all long roots and so the  highest long root is the highest  weight of $\u_t.$
\medskip

\textbf{\emph{\underline{Type $C:$}}} In this case, we have
$\a_{sh}^1=\ep_{n-1}+\ep_n,$ $\a_{sh}^2=\ep_{\ell-1}+\ep_\ell,$
$\a_{lg}^1=2\ep_{n}$ and $\a_{lg}^2=2\ep_\ell.$ Setting
$\b:=\ep_{\ell-1}-\ep_\ell,$ we get $\a_{sh}^2+\b\in(R_1)_{lg}.$
Now using the same argument as in Type $B,$ we are done.

\medskip

\textbf{\emph{\underline{Type $BC:$}}}  In this case, we have
$\a_{sh}^1=\ep_n,$ $\a_{sh}^2=\ep_\ell,$
$\a_{lg}^1=\ep_{n-1}+\ep_n,$  $\a_{lg}^2=\ep_{\ell-1}+\ep_\ell,$
$\a_{ex}^1=2\ep_{n}$ and $\a_{ex}^2=2\ep_\ell.$ One can also easily
see that  $$\{\gamma+\b\mid \gamma\in(R_1)_{sh},\b\in(R_1)_{lg}\cup
(R_{1})_{ex}\}\cap R_1\sub(R_1)_{sh}.$$ This together with
(\ref{tavakol1}), (\ref{gen3}) and the fact that $\D_1\sub (R_1)_{ex}\cup(R_1)_{lg}$ proves the claim stated in  Step 1 in the case that
$\a=\a_{sh}^2.$ Now suppose that $\a=\a_{lg}^2.$ Setting
$\b:=\ep_{\ell-1}-\ep_\ell,$ we get  that $\a+\b$ is an extra long root.
Now   we are done using the same argument as in Type $B.$ Next
suppose $\a=\a_{ex}^2,$ then by (\ref{gen3}), $u_t$ is a weight
vector of weight $\a$ which is an extra long root. Therefore  any
extra long root is a weight for $\u_t$  and so the highest weight of
$\u_t$ is $\a_{ex}^1.$ This completes the proof of Step 1.
\medskip

\noindent{\bf Step 2:} We first note that depending on the type of $R_2,$ $\a$ is one of $\ep_\ell,2\ep_\ell,\ep_\ell+\ep_{\ell-1}$ or $\ep_{\ell+1}-\ep_{1}.$
If $\a=\ep_{\ell}+\ep_{\ell-1},$ then either $R_2$ is of type $B$ or $D$ and  by Step 1, the set of nonzero  weights of $\u$ coincides with $R^\times,$ or $R$ is of type $C$ or $BC$ and the set of nonzero weights of $\u$ coincides with $(R_{sdiv})_{sh}.$  In both cases,  using induction on $r\in\bbbn\setminus\{0\},$ one can see that
\begin{equation}\label{rev1}
\parbox{4in}{
 if  $1\leq m_1,\ldots,m_r\leq n$  and for each $1\leq p\leq r,$ $\a_{m_p}+\cdots+\a_{m_1}+\a$ is a weight  for $\u,$ then     $\{m_1,\ldots,m_r\}\sub\{\ell,\ldots,n\}$ and  $\a_{m_r}+\cdots+\a_{m_1}+\a=\ep_{q}+\ep_{q'}$ for some $\ell-1\leq q\neq q'\leq n.$
}\end{equation}
Also if  $\a=\ep_{\ell+1}-\ep_1,\ep_\ell,2\ep_{\ell},$ one  can see that
\begin{equation}\label{rev2}\parbox{3.5in}{
if $r$ is a positive integer and  $1\leq m_1,\ldots,m_r\leq n$ are such that  for each $1\leq p\leq r,$ $\a_{m_p}+\cdots+\a_{m_1}+\a$ is a weight  for $\u,$ then  $\{m_1,\ldots,m_r\}\sub\{\ell+1,\ldots,n\}.$ }\end{equation}

Now suppose that    $0\neq u\in \u_{\a}, $ we shall  show that $u$ is a
scalar multiple of $v.$ Since $u\in \u,$ by (\ref{tavakol1}), $u$ is
written as a linear combination of weight vectors  of the form
$f_{i_t}\cdot \cdots\cdot  f_{i_1},e_{j_s}\cdot \cdots\cdot  e_{j_1}\cdot v,$
$t,s\in\bbbn,\;1\leq i_1,\ldots,i_t,j_1,\ldots,j_s\leq n.$ So
without loss of generality, we  suppose
$$u=f_{i_t}\cdot\cdots\cdot f_{i_1}\cdot e_{j_s}\cdot \cdots\cdot e_{j_1}\cdot v$$ where
$t,s\in\bbbn,$ and $1\leq i_1,\ldots, i_t,j_1,\ldots,j_s\leq n.$
Since $u$ is of weight $\a,$ we get that
$\a+\a_{j_1}+\cdots+\a_{j_s}-\a_{i_1}-\cdots-\a_{i_t}=\a.$ This
implies that
\begin{equation}
\label{coin} s=t\andd
(j_1,\ldots,j_s)=(\sg(i_1),\ldots,\sg(i_t))
\end{equation}
for a permutation $\sg$ of $\{i_1,\ldots,i_t\}.$
We note that  $\a$ is an element of $R_{2}^+$ and so it  is written
as a linear combination of $\{\a_i\mid 1\leq i\leq
\ell\}$ with nonnegative rational  coefficients not all equal to zero. Now since  $\{\a_i\mid 1\leq
i\leq n\}$ is a base of $(R_{1})_{sdiv},$ $\a-\a_{j}$ ($\ell+1\leq j\leq n$)
is not a root of $R_1$ and so it is not an element of $\Lam_1.$
Therefore
\begin{equation}\label{gen5}f_{j}\cdot v=0,\;\;\ell+1\leq j\leq n.
\end{equation} Now this
 implies that
\begin{equation}
\label{eq1}\begin{array}{c} f_{j}\cdot e_{j}\cdot v=e_j\cdot
f_j\cdot v-h_{j}\cdot v=0-h_{j}\cdot v\in\bbbf v,\\
(\ell+1 \leq j\leq n).\end{array}
\end{equation}
We also note that as $v$ is  a highest vector of $\gg_2$-module
$\w,$ $e_j\cdot v=0$ for $ 1\leq j\leq \ell.$ Therefore one gets that
\begin{equation}\label{js}
\hbox{$j_1\in\{\ell+1,\ldots, n\}$ provided that  $s\neq0.$ }
\end{equation}

Now we are ready to prove that  $u$ is a scalar multiple of $v.$
If $s=0,$ there is nothing to prove. So we suppose  $s\geq1,$ and
use induction on $s$ to prove. If $s=1,$ we get the result
appealing (\ref{coin}), (\ref{js}) and (\ref{eq1}). Now suppose
$s>1.$ If $\a=\ep_{\ell}+\ep_{\ell-1},$ then using  (\ref{rev1}) together with  (\ref{coin}), we get that $i_1\in\{\ell,\ldots,n\}.$ This together with (\ref{gen5}) and  the fact that $2\ep_{\ell-1}$ is not a weight for $\u$ implies that
$f_{i_1}\cdot v=0.$ Next  take $1\leq k_1<\ldots<
k_r\leq s$ to be the only indices with
$j_{k_1}=\cdots=j_{k_r}=i_1$ and  use  Lemma \ref{gen-fact} to get
$$f_{i_t}\cdot\cdots\cdot f_{i_1}\cdot e_{j_s}\cdot\cdots\cdot e_{j_1}\cdot v\in\sum_{q=1}^r\bbbf f_{i_t}\cdot\cdots\cdot f_{i_2}\cdot e_{j_s}\cdot\cdots e_{j_{k_r}}\cdot\cdots \cdot \widehat{e_{j_{k_q}}}\cdot\cdots\cdot e_{j_{k_1}}\cdot\cdots \cdot e_{j_1}\cdot v.$$
Now induction hypothesis completes the proof of this step in the case that $\a=\ep_{\ell-1}+\ep_\ell.$ Next suppose $\a\in\{\ep_\ell,2\ep_{\ell},\ep_{\ell+1}-\ep_{1}\},$ then using  (\ref{rev2})  together with  (\ref{coin}), (\ref{gen5}), Lemma \ref{gen-fact}  and the same argument as before, we  get  the result.
\medskip

{\bf Step 3:} It is immediate using Step 2 together with
(\ref{gen3}).\qed

\section{Root graded Lie algebras} The structure of  Lie algebras graded by an irreducible  finite root system have been studied in \cite{BM}, \cite{BZ}, \cite{ABG1}, \cite{N}, \cite{ABG2} and \cite{BS}. A Lie algebra $\LL$ graded by an irreducible  finite root system $R$ contains a finite dimensional  split simple  Lie algebra $\gg$ and with respect to a splitting Cartan subalgebra, it  is equipped with a weight space decomposition whose  set of weights is contained in $R.$ This feature  allows us to decompose $\LL$ into finite dimensional  irreducible $\gg-$submodules whose set of nonzero weights is $R_{sh},$ $R_{sdiv}^\times$ or $(R_{sdiv})_{sh}.$ Collecting the components  of the same highest weight results in  the  decomposition
\begin{equation}\label{*}\LL=(\gg\ot \aa)\op(\ss\ot\bb)\op(\v\ot \cc)\op \dd\end{equation}  in which $\dd$ is a trivial submodule of $\LL$ and $\ss$ (resp. $\v$) is  the finite dimensional irreducible $\gg-$module whose  set of nonzero weights is  $(R_{sdiv})_{sh}$  (resp.  $R_{sh}$).   The Lie algebraic structure on $\LL$ induces an algebraic structure on $\fb:=\aa\op\bb\op\cc$ which we refer to  as the {\it coordinate algebra} of $\LL.$ The Lie bracket on $\LL$ can be rewritten using the ingredients involved in describing the product defined on  the algebra $\fb.$ In this section, we have a compromising  view on the coordinate algebras of root graded  subalgebras of $\LL.$ We devote this section to two subsections. In the first subsection, we illustrate the structure of a specific Lie algebra which we shall frequently   use in   the sequel of the paper.
In the second  subsection,  we consider a Lie algebra $\LL$ graded by an irreducible  finite
root system $R$ and for a full irreducible subsystem $S$ of $R$ which is  of the same type as $R,$ we take   $\LL^{^S}$ to be the Lie subalgebra of $\LL$ generated by homogeneous spaces in correspondence to $S^\times.$ We show that the coordinate algebra of the Lie subalgebra  $\LL^{^S},$ which is an algebra graded by $S,$ does not depend on  $S.$ In fact, we prove that the coordinate algebra of $\LL^{^S}$ coincides with
the coordinate algebra of $\LL.$ Moreover, we describe the Lie bracket on $\LL$ in terms of  the ingredients involved in describing the Lie bracket on $\LL^{^S}$ with respect to its coordinate algebra. Our method is based on a  type-by-type approach. Since the proofs for different  types are quite similar, we go  through the proofs in details if $R$ is of type $BC$ and  for other types, we just report the results and leave the proofs  to the readers.

\subsection{A specific Lie algebra}\label{subsect2-1}
By a {\it star algebra} $(\fa,*),$ we mean an algebra $\fa$ together
with   a self-inverting antiautomorphism $*$ which is referred  to as an {\it involution}. We call a quadruple $(\fa,*,\cc,f)$ a {\it coordinate  quadruple} if one of the following holds:
\begin{itemize}
\item  (Type $A$) $\fa$ is a unital associative algebra, $*=id_\fa,$ $\cc=\{0\}$ and $f:\cc\times\cc\longrightarrow \fa$ is the  zero map.
\smallskip

\item    (Type $B$) $\fa=\aa\op\bb$  where $\aa$ is a unital commutative associative algebra and $\bb$ is a unital associative  $\aa-$module equipped with a symmetric bilinear form and $\fa$ is the corresponding Clifford  Jordan algebra, $*$ is a linear transformation fixing the elements of $\aa$ and skew fixing the elements of $\bb,$ $\cc=\{0\}$ and $f:\cc\times\cc\longrightarrow \fa$ is the zero map.
\smallskip

\item  (Type $C$) $\fa$ is a unital associative algebra, $*$ is an involution, $\cc=\{0\}$ and $f:\cc\times\cc\longrightarrow \fa$ is the zero map.
\smallskip

\item  (Type $D$) $\fa$ is a unital  commutative associative algebra $*=id_\fa,$ $\cc=\{0\}$ and $f:\cc\times\cc\longrightarrow \fa$ is the zero map.
\smallskip

\item  (Type $BC$) $\fa$ is a unital associative algebra, $*$ is an involution, $\cc$ is a unital associative   $\fa-$module and $f:\cc\times\cc\longrightarrow \fa$ is a skew-hermitian form.

\end{itemize}

Suppose that $(\fa,*,\cc,f)$ is a coordinate quadruple. Denote by  $\aa$ and $\bb,$ the fixed and the skew fixed points of  $\fa$ under $*,$ respectively. 
 Set
$\fb:=\fb(\fa,*,\cc,f):=\fa\op\cc$ and define
\begin{equation}\label{probinbc-n}
\begin{array}{c}\cdot:\fb\times\fb\longrightarrow \fb\\
(\a_1+c_1,\a_2+c_2)\mapsto(\a_1\cdot \a_2)+f(c_1,c_2)+\a_1\cdot
c_2+\a_2^*\cdot c_1,
\end{array}
\end{equation}
for $\a_1,\a_2\in\fa$ and $c_1,c_2\in\cc.$
Also for $\b,\b'\in\fb,$
set
\begin{equation}\label{end9-n}
\b\circ\b':=\b\cdot\b'+\b'\cdot\a\andd [\b,\b']:=\b\cdot \b'-\b'\cdot\b
\end{equation}
and for $c,c'\in\cc,$ define
\begin{equation}\label{diamond-heart-n}
\begin{array}{c}
\begin{array}{ll}
\diamond:\cc\times\cc\longrightarrow\aa,&
(c,c')\mapsto\frac{f(c,c')-f(c',c)}{2};\;c,c'\in\cc,
\end{array}\vspace{3mm}\\
\begin{array}{ll}
\heart:\cc\times\cc\longrightarrow\bb,&
(c,c')\mapsto\frac{f(c,c')+f(c',c)}{2};\;c,c'\in\cc.
\end{array}
\end{array}
\end{equation}

Now suppose that $\ell$ is a positive integer and for $\a,\a'\in\fa$ and $c,c'\in\cc,$ consider the following
endomorphisms
\begin{equation}\label{derivbc}
\begin{array}{l}
d_{\a,\a'}^{\ell,\fb}:\fb\longrightarrow\fb,\\
\b\mapsto\left\{\begin{array}{ll}
\frac{1}{\ell+1}[[\a,\a'],\b]& X=A_\ell,\; \b\in \fb,\vspace{2mm}\\
\a'(\a\b)-\a(\a'\b)&X=B_\ell,\;\b\in\fb,\vspace{2mm}\\
 \frac{1}{4\ell}[[\a,\a']+[\a^{*},\a'^{*}],\b]&X=C_\ell,BC_\ell,\;\;\b\in\fa,\vspace{2mm}\\
\frac{1}{4\ell}([\a,\a']+[\a^{*},\a'^{*}])\cdot
\b&X=C_\ell,BC_\ell,\;\;\b\in\cc,\vspace{2mm}\\
0& X=D_\ell,\; \b\in\fb,\end{array}\right.\\d^{\ell,\fb}_{c,c'}:\fb\longrightarrow\fb,\\
\b\mapsto\left\{\begin{array}{ll}\frac{-1}{2\ell}[c\heart c',\b]&\; X=BC_\ell,\;\b\in\fa,\vspace{2mm}\\
\frac{-1}{2\ell}(c\heart c')\cdot\b-\frac{1}{2}(f(\b,c')\cdot
c+f(\b,c)\cdot c')&X=BC_\ell,\;\b\in\cc,\vspace{2mm}\\
0&\hbox{otherwise},\vspace{2mm}\end{array}\right.\\
d^{\ell,\fb}_{\a,c}:=d^{\ell,\fb}_{c,\a}:=0,\vspace{2mm}\\
d^{\ell,\fb}_{\a+c,\a'+c'}:=d^{\ell,\fb}_{\a,\a'}+d^{\ell,\fb}_{c,c'}.
\end{array}
\end{equation}
One can see that for $\b,\b'\in\fb,$ $d^{\ell,\fb}_{\b,\b'}\in Der(\fb).$
Next take  $K$ to be  a subspace of $\fb\ot \fb$ spanned by
$$\begin{array}{l}
\a\ot c,\;\;c\ot\a,\;\;a\ot b,\\
\a\ot\a'+\a'\ot\a,\;\;c\ot c'-c'\ot c,\\
(\a\cdot \a')\ot\a''+(\a''\cdot\a)\ot\a'+(\a'\cdot\a'')\ot\a,\\
f(c,c')\ot\a+( \a^*\cdot c')\ot c-(\a\cdot c)\ot c'
\end{array}$$
for $\a,\a',\a''\in\fa,$ $a\in\aa,$ $b\in\bb,$ and $c,c'\in\cc.$
Then $(\fb\ot\fb)/K$ is a Lie algebra under the Lie bracket
\begin{equation}\label{last4}
[(\b_1\ot \b_2)+K,(\b'_1\ot \b'_2)+K]_\ell=((d^{\ell,\fb}_{\b_1,\b_2}(\b'_1)\ot \b'_2)+K)+(\b'_1\ot d^{\ell,\fb}_{\b_1,\b_2}
(\b'_2))+K)\end{equation}
for $\b_1,\b_2,\b'_1,\b'_2\in\fb$ (see \cite[Proposition 5.23]{ABG2} and \cite{ABG1}). We denote this Lie algebra by  $\{\fb,\fb\}_\ell$  (or $\{\fb,\fb\}$ if there is no confusion) and for
$\b_1,\b_2\in\fb,$ we  denote $(\b_1\ot\b_2)+K$ by   $\{\b_1,\b_2\}_\ell$ (or $\{\b_1,\b_2\}$  if there is no confusion). We recall the {\it full skew-dihedral homology group}  $$FH(\fb):=\{\sum_{i=1}^n\{\b_i,\b'_i\}_\ell\in\{\fb,\fb\}_\ell\mid \sum_{i=1}^nd^{\ell,\fb}_{\b_i,\b'_i}=0\}$$  of $\fb$  (with respect to $\ell$) from \cite{ABG2} and \cite{ABG1} and  note that it is a subset of the center of $\{\fb,\fb\}_\ell. $
For $\b_1=a_1+b_1+c_1\in\fb$ and $\b_2=a_2+b_2+c_2\in\fb$ with $a_1,a_2\in\aa,$ $b_1,b_2\in\bb$ and $c_1,c_2\in\cc,$ set \begin{equation}\label{beta*}\b_{_{\b_1,\b_2}}^*:=[a_1,a_2]+[b_1,b_2]-c_1\heart c_2.\end{equation}
We say a subset   $\kk$ of the full skew-dihedral homology group of $\fb$ satisfies the ``{\it uniform property on $\fb$}" if  for  $\b_1,\b'_1,\ldots,\b_n,\b'_n\in\fb,$ $\sum_{i=1}^n\{\b_i,\b'_i\}_\ell\in\kk$   implies that  $\sum_{i=1}^n\b^*_{\b_i,\b'_i}=0.$
\begin{rem}\label{rem3}{\rm Suppose that $\ell,\ell'$ are two positive integers. If $\kk$ is a subset of the full skew-dihedral homology group of $\fb(\fa,*,\cc,f)$ with respect to $\ell$ satisfying the  uniform property on $\fb(\fa,*,\cc,f)$, it  is a subset of the full skew-dihedral homology group of $\fb(\fa,*,\cc,f)$ with respect to $\ell'$ satisfying the  uniform property on $\fb(\fa,*,\cc,f).$ In other words, the uniform property on $\fb(\fa,*,\cc,f)$ dose  not depend on $\ell.$}
\end{rem}


%

\subsection{Lie algebras graded by a finite root system}\label{subsect2-2}
In this  work, we study root graded Lie algebras in the following sense:

\begin{Definition}\label{root-graded}{\rm
Suppose that $R$ is  an irreducible locally finite root system. We
say a Lie algebra $\LL$ is an {\it $R-$graded Lie algebra with graded
pair $(\gg,\hh)$} if the followings are satisfied:

i) $\gg$ is  a locally finite split simple Lie subalgebra of $\LL$ with
splitting Cartan subalgebra $\hh$  and corresponding root system
$R_{sdiv}.$

ii) $\LL$ has  a weight space decomposition $\LL=\op_{\a\in
R}\LL_\a$ with respect to $\hh$ via the adjoint representation.

iii) $\LL_0=\sum_{\a\in \rcross}[\LL_\a,\LL_{-\a}].$}
\end{Definition}

The following lemma easily  follows from Lemma \ref{final3}.
\begin{Lemma}\label{final1}
Suppose that $R$ is an irreducible  locally finite root system and
$\LL$ is a Lie algebra graded by $R$ with grading pair
$(\gg,\hh).$ Let $S$ be an irreducible full subsystem of $R$ and
set
\begin{equation}\label{final2}
\begin{array}{l}
\displaystyle{\LL^{^S}:=\sum_{\a\in S^\times}\LL_\a\op\sum_{\a\in
S^\times}[\LL_\a,\LL_{-\a}],}\\
\displaystyle{\gg^{^S}:=\sum_{\a\in
S_{sdiv}^\times}\gg_\a\op\sum_{\a\in
S_{sdiv}^\times}[\gg_\a,\gg_{-\a}].}
\end{array}
\end{equation}
Then $\LL^{^S}$ is an $S-$graded Lie subalgebra of $\LL$ with
grading pair $(\gg^{^S},\hh^{^S}:=\hh\cap\gg^{^S}).$
\end{Lemma}

Before going through the main body of this subsection, we  want to fix a notation. If  $A$ is a subspace of a vector space $V_1$ and $B$ is a subspace of a vector space  $V_2,$ by a conventional notation, we take $A\dot\ot B$ to be  the vector subspace
of $V_1\ot V_2$  spanned by $a\ot b$ for $a\in A$ and $b\in B.$

\subsubsection{\textbf{Type $BC$}}\label{subsub1}
Suppose that $I$ is a nonempty index set  of cardinality
$m_n:=n>3$ and  $I_0$ is a nonempty subset of $I$ of cardinality
$m_\ell:=\ell>3.$ Take $\v:=\v^n$ to be a vector space with a basis
$\{v_i\mid i\in I\cup\bar I\}$  equipped with a nondegenerate
symmetric bilinear form  $\fm$ as in (\ref{form-c}). Set $\gg^n:=\mathfrak{sp}(I)$ and take $\ss:=\ss^n$ to be as in
(\ref{module-s-c}). Consider (\ref{simple-c}), (\ref{simple-c-alg}) and set
$$\v^\ell:=\v_{_{I_0}},\;\;\gg^\ell:=\gg_{_{I_0}},\;\;\ss^\ell:=\ss_{_{I_0}}.$$

We  take $Id_{_\v}$ to be the identity map on $\v$ and define the linear endomorphism  $Id_{_{\v^\ell}}$ on $\v$ by
$$\begin{array}{c}Id_{_{\v^{\ell}}}:\v\longrightarrow \v\\
v_i\mapsto v_i,\;v_{\bar i}\mapsto v_{\bar i},\;v_{j}\mapsto 0,\;v_{\bar j}\mapsto 0;\;\;
(i\in I_0,\; j\in I\setminus I_0).
\end{array}$$
Also for $\lam=\ell,n$ and $x,y\in\gg^\lam\op\ss^\lam,$  set
\begin{equation}
\label{end7}
x\circ_\lam y:=xy+yx-(1/m_\lam)tr(xy)Id_{_{\v^\lam}}.
\end{equation}

Next for $u,v\in\v,$ define
\begin{equation}\label{u,v}\begin{array}{l}
\;[u, v]
:\v\longrightarrow\v;\;w\mapsto\frac{1}{2}((v,w)u+(w,u)v)+\frac{1}{2\ell
}(u,v)Id_{_{\v^\ell}}(w);\;\;w\in\v,\vspace{2mm}\\
\;u\circ v:\v\longrightarrow\v;\;w\mapsto\frac{1}{2}((v,w)u+(u,w)v);\;\;w\in\v,\vspace{2mm}\\
\;[u, v]_n
:\v\longrightarrow\v;\;w\mapsto\frac{1}{2}((v,w)u+(w,u)v)+\frac{1}{2n
}(u,v)Id_{_\v}(w);\;\;w\in\v.
\end{array}
\end{equation}

One can easily see that  up to isomorphism  $$\begin{array}{ll}\gg^\ell=\hbox{span}\{u\circ  v\mid u,v\in\v^\ell\},&
\ss^\ell=\hbox{span}\{[u,v]\mid u,v\in\v^\ell\},\\\gg^n=\hbox{span}\{u\circ  v\mid u,v\in\v^n\},&
\ss^n=\hbox{span}\{[u,v]_n\mid u,v\in\v^n\}.\end{array}$$

Suppose that $R$ is an irreducible finite   root system of type $BC_I$ and $S$ is the  irreducible full subsystem of $R$ of
type $BC_{I_0}.$ Suppose that     $\LL$ is an $R-$graded Lie algebra with grading pair $(\fg,\fh)$ and take $\LL^{^S},$ $\fg^{^S}$ and $\fh^{^S}$ to be  as in Lemma  \ref{final1}. In order to simplify the  using of the  notations, we set
\begin{equation}\label{general}
\LL^n:=\LL,\;\; \LL^{\ell}:=\LL^{^S},\;\;[u,
v]_\ell:=[u,v];\;\;(u,v\in\v).
\end{equation}

One knows that as a $\fg^{^S}-$module, $\LL^\ell$  can be decomposed into finite dimensional  irreducible $\fg^{^S}-$submodules,
each of which is a  finite dimensional  irreducible $\fg^{^S}-$module with highest weight contained in $S.$ Take
\begin{equation}\label{last3-fin}\LL^\ell=\LL^{^S}=\bigoplus_{i\in
\i_0}\fg_i\op\bigoplus_{j\in\jj_0}\mathfrak{s}_j\op\bigoplus_{t\in
\T_0}V_t\op E\end{equation} to be the decomposition of $\LL^\ell$
into finite dimensional irreducible $\fg^{^S}-$modules in which
$\i_0,\jj_0,\T_0$ are (possibly empty) index sets and for
$i\in\i_0,$ $j\in \jj_0,$ and $t\in\T_0,$ $\fg_i$ is isomorphic to
$\fg^{^S}(\simeq\gg^\ell),$  $\sfs_j$ is isomorphic to $\ss^\ell,$
$V_t$ is isomorphic to $\v^\ell$ and $E$ is a trivial
$\fg^{^S}-$submodule.
\begin{Lemma}\label{divide1}
Use the notation as in the text and consider $\LL=\LL^n$ as a
$\fg-$module. Then there exist index sets $\i,$ $\jj,$ $\T$  with
$\i_0\sub\i,\jj_0\sub\jj,\T_0\sub\T,$  and a class $\{\dd_n,\gg_i,\ss_j,\v_t\mid i\in \i,j\in \jj,t\in \T\}$ of  finite dimensional   $\fg-$submodules of $\LL$ such that

\begin{itemize}
\item   $\dd_n$ is a trivial $\fg-$module,  $\gg_i$  is  isomorphic to $\fg,$ $\ss_j$
 is  isomorphic to $\ss,$ and $\v_t$
is isomorphic to $\v,$  for $i\in\i,j\in \jj,$  $t\in \T,$  \\
\item $\fg_i\sub\gg_i,$   $\sfs_j\sub\ss_j,$
$V_t\sub\v_t$  ($i\in\i_0,$ $j\in\jj_0,$ $t\in\T_0$),\\ \item
$\LL^n=\LL=\bigoplus_{i\in\i}\gg_i\op\bigoplus_{j\in\jj}\ss_j\op\bigoplus_{t\in\T}\v_t\op\dd_n.$
\end{itemize}We make a convention that we refer to \begin{equation}\label{last6}(\ii,\jj,\T,\{\fg_i\},\{\gg_i\},\{\sfs_j\},\{\ss_j\},\{V_t\},\{\v_t\},E,\dd_n)\end{equation} as  an $(R,S)-$datum for the pair $(\LL^n,\LL^\ell).$

\end{Lemma}

\pf For $i\in \ii_0,$ by Proposition \ref{gen1}, the $\gg-$submodule
$\gg_i$ of $\LL $ generated by $\fg_i$ is a finite dimensional
$\gg-$module isomorphic to $\gg.$ For $\a\in
S_{sdiv}^\times$ and $0\neq x\in(\fg_i)_\a\sub\LL_\a,$ we have
$x\in\LL_\a\cap\fg_i\sub\LL_\a\cap\gg_i=(\gg_i)_\a.$
Now as $\dim(\gg_i)_\a=\dim(\fg_i)_\a=1,$ we get
\begin{equation}\label{final4}(\fg_i)_\a=(\gg_i)_\a;\;\;\;\a\in S^\times_{sdiv}.\end{equation}
 \underline{\noindent\textbf{Claim 1.}} The sum $\sum_{i\in \mathcal{I}_0
}{\gg} _i$ is a direct sum: Suppose that
$i_0,i_1,\ldots,i_n$ are distinct elements of $\ii_0$ and $0\neq x \in\gg_{i_0}
\cap\sum_{t=1}^n\gg_{i_t}.$ Then as $\gg_{i_0} $ is an irreducible $\gg
-$module, we get that
$$\gg_{i_0} \sub\sum_{t=1}^n\gg_{i_t}.$$  This together with (\ref{final4}) implies that for  $\a\in S_{sdiv}^\times\sub R_{sdiv}^\times,$
$$(\fg_{i_0})_\a=(\gg_{i_0})_\a\sub\sum_{t=1}^n(\gg_{i_t})_\a=\sum_{t=1}^n(\fg_{i_t})_\a$$ which
contradicts the fact that $\sum_{i\in \mathcal{I}_0}{\fg}_i$ is
direct. This completes the proof of  Claim 1.

\smallskip

Now for $j\in\jj_0$
and $t\in\mathcal{T}_0,$ take $\ss _j$ and $\v _t$ to be  the
finite dimensional irreducible $\gg-$submodules of $\LL $
generated  by  $\sfs_j$ and $V_t$ respectively. Using the same
argument as above, one can see that the summations $\sum_{j\in\jj_0}\ss _j$ and
$\sum_{t\in\mathcal{T}_0}\v _t$ are direct. Set
$$\gg(n):=\op_{i\in \mathcal{I}_0
}{\gg} _i,\;\;\ss(n):=\op_{j\in\jj_0}\ss _j,\;\;\v(n):=\op_{t\in\mathcal{T}_0}\v _t.$$
We note  that $\gg(n)$ (resp. $\ss(n)$ and $\v(n)$) is a  $\gg-$submodule of $\LL$ whose set of weights is  $R_{sdiv}$ (resp. $R_{lg}\cup\{0\}$ and $R_{sh}$).
\medskip

\underline{\noindent\textbf{Claim 2.}} For $\a\in R_{lg},$ $x\in \gg(n)_\a,$ and
$y\in \ss(n)_\a,$ we have  $x+y=0$ if and only if  $x=y=0:$ Suppose that $x+y=0.$ Since
$x\in\gg(n)_\a=\sum_{i\in\ii_0}(\gg_i )_\a,$ we get
$x=\sum_{i\in\ii_0}x_i$ with finitely many nonzero terms
$x_i\in(\gg_i)_\a,$ for $i\in\ii_0.$ Similarly
$y=\sum_{j\in\jj_0}y_j$ with finitely many nonzero terms
$y_j\in(\ss_j)_\a,$ for $j\in\jj_0.$ Now we recall that $\ell,n>3$ and  $R,S$ are  root systems of type $BC_n$ and $BC_\ell$ respectively. This allows us to pick $\b_1,\b_2\in R_{lg}$ such that
$\b:=\a+\b_1+\b_2\in S_{lg}$ and $\a+\b_1\in R_{lg}.$ Fix $a_1\in\gg_{\b_1}$ and $a_2\in\gg_{\b_2}.$
Using Lemma \ref{elementary}, we get that
\begin{equation}\label{eq2}
\parbox{3.6in}{ if
$x_i$ $(i\in\ii_0)$ is nonzero, then $ a_2\cdot a_1\cdot x_i$ is a nonzero element of $(\gg_i
)_\b=(\fg_i)_\b$ and similarly  if $j\in\jj_0$ and  $y_j\neq0,$
$ a_2\cdot a_1\cdot y_j$ is  a nonzero element of $(\ss_j )_\b=(\sfs_j)_\b.$ }\end{equation}Now
since $x+y=0,$ we get that
$\sum_{i\in\ii_0}x_i=-\sum_{j\in\jj_0}y_j$ which in turn implies that
$\sum_{i\in\ii_0}a_2\cdot a_1\cdot x_i=-\sum_{j\in\jj_0} a_2\cdot
a_1\cdot y_j.$ But the right hand side is an element of
$\op_{j\in\jj_0}(\sfs_j)_\b$ and the left hand side is an element
of $\op_{i\in\ii_0}(\fg_i)_\b.$
Therefore $a_2\cdot a_1\cdot x_i=0$ and $a_2\cdot
a_1\cdot y_j=0$ for $i\in\ii_0$ and $j\in\jj_0.$ This together with (\ref{eq2}) implies
that for $i\in \ii_0$ and $j\in\jj_0,$ $x_i=0$ and $y_j=0.$   This completes the proof of Claim 2.

\medskip

\underline{\noindent\textbf{Claim 3.}} For $x\in\gg(n)_0$ and $y\in\ss(n)_0,$
$x+y=0$ if and only if $x=y=0:$ Suppose that $x+y=0$ and  $x\neq0.$ Since
$\gg(n)_0= \sum_{i\in \ii_0}(\gg_i)_0,$ we have
$x=\sum_{i\in\ii_0}x_i$ with finitely many nonzero terms  $x_i\in
(\gg_i)_0,$ $i\in\ii_0.$ Fix  $t\in\ii_0$ such that
$x_t\neq 0.$ Since  $x_t$ is a nonzero element of the
irreducible nontrivial $\fg -$module $\gg_t ,$ there is
 $\a\in R^\times$ and $0\neq a\in\fg_\a$ such that $a\cdot x_t\neq
 0.$ We  note  that $x\in\gg(n)_0$ and
$y\in\ss(n)_0,$  therefore we have $a\cdot x\in\gg(n)_\a$
and $a\cdot y\in\ss(n)_\a.$ Now as $0=a\cdot x +a\cdot
y,$ Claim 2 together with the fact that the set of weights of
$\ss(n)$ is $R_{lg}\cup\{0\}$ implies that $a\cdot x=0$ and
$a\cdot y=0.$ So $\sum_{i\in\ii_0}a\cdot x_i=0.$ But by Claim 1,
$\sum_{i\in\ii_0}\gg_i$ is a direct sum, so $a\cdot x_t=0$ which
is a contradiction. Therefore $x=0$ and so $y=0$ as well. This completes the proof of Claim 3.

\smallskip

\underline{\textbf{Claim 4.}}  The sum $\gg(n)+\ss(n)+\v(n)$ is a direct sum:
Suppose that $x\in\gg(n),$ $y\in\ss(n)$  and $z\in\v(n)$ are such
that $x+y+z=0.$ We have  $x=\displaystyle{\sum_{\a\in
R_{sdiv}}x_\a}$ with $x_\a\in\gg(n)_\a\sub\LL_\a$ for  $\a\in
R_{sdiv},$ $y=\displaystyle{\sum_{\a\in R_{lg}\cup\{0\}}y_\a}$
with $y_\a\in\ss(n)_\a\sub\LL_\a$ for  $\a\in R_{lg}\cup\{0\},$
and $z=\displaystyle{\sum_{\a\in R_{sh}}z_\a}$ with
$z_\a\in\v(n)_\a\sub\LL_\a$ for  $\a\in R_{sh}.$ Therefore one
gets that
$$
\begin{array}{c}
x_0+y_0=0,\;
z_\a=0,\;x_\b+y_\b=0, \;x_\gamma=0;\\( \a\in R_{sh},\;\b\in R_{lg},\;\gamma\in R_{ex}).
\end{array}$$
Now using Claims 2,3, we are done

\medskip
To complete the proof, we note that  as a $\fg-$module, $\LL $ can be decomposed into  finite dimensional irreducible $\fg-$submodules with  the set of weights contained in $R .$ Now as  $\bigoplus_{i\in\ii_0}\gg
_i\op\bigoplus_{j\in\jj_0}\ss _j\op\bigoplus_{t\in\mathcal{T}_0}\v
_t$ is a submodule of $\LL ,$   one can find index sets $\ii ,\jj
,\mathcal{T}$ with
$$\ii_0\sub\ii ,\;\jj_0\sub\jj ,\;\mathcal{T}_0\sub\mathcal{T} $$
and a class $\{\dd_n,\gg_i,\ss_j,\v_t\mid i\in\ii\setminus\ii_0,j\in\jj\setminus\jj_0,t\in\T\setminus\T_0\}$ of finite dimensional  $\fg-$submodules such that $\dd_n$ is a trivial $\fg-$module, $\gg_i $
is  isomorphic to $\gg$ $(i\in\ii\setminus\ii_0),$ $\ss_j$ is  isomorphic to $\ss$ ($j\in\jj\setminus \jj_0$),  $\v _t$
is isomorphic to $\v $ $(t\in\mathcal{T} \setminus\mathcal{T}_0)$ and
\begin{eqnarray*}
\LL &=&(\bigoplus_{i\in\ii_0}\gg _i\op\bigoplus_{j\in\jj_0}\ss
_j\op\bigoplus_{t\in\mathcal{T}_0}\v _t) \op(\bigoplus_{i\in\ii\setminus\ii_0}\gg _i\op\bigoplus_{j\in\jj\setminus\jj_0}\ss
_j\op\bigoplus_{t\in\mathcal{T}\setminus\T_0}\v _t \op\dd_n)\\&=&\bigoplus_{i\in\ii}\gg _i\op\bigoplus_{j\in\jj}\ss
_j\op\bigoplus_{t\in\mathcal{T}}\v _t \op\dd_n.
\end{eqnarray*}
This completes the proof.\qed

\bigskip

From now on, we use the data appeared in the $(R,S)-$datum (\ref{last6}). We take $\aa_n $ to be a vector space with a basis $\{a_i\mid i\in
\ii \},$   $\bb_n $ to be a vector space with a basis $\{b_j\mid
t\in\jj \}$ and $\cc_n $ to be a vector space with a basis
$\{c_t\mid t\in \mathcal{T}\}.$  Then as a $\gg^n-$module, $\LL $
can be identified with \begin{equation}\label{last1}(\gg^n \ot \aa_n )\op(\ss^n \ot \bb_n )\op(\v^n
\ot\cc_n )\op\dd_n .\end{equation} Take
\begin{equation}\label{identification}\varphi:\LL\longrightarrow (\gg^n \ot \aa_n )\op(\ss^n \ot \bb_n )\op(\v^n
\ot\cc_n )\op\dd_n \end{equation} to be the canonical  identification. Next define $\aa_\ell$ to be the vector subspace
of $\aa_n $ spanned by $\{a_i\mid i\in\ii_0\},$ $\bb_\ell$ to be the
vector subspace of $\bb_n $ spanned by $\{b_j\mid
j\in\mathcal{J}_0\}$ and $\cc_\ell$ to be  the  vector subspace of $\cc_n $
spanned by $\{c_t\mid t\in \mathcal{T}_0\}.$ Then it follows from (\ref{last3-fin}) that as a $\gg^\ell-$module,
$\LL^\ell=\LL^{^S}$ can be identified with
\begin{equation}\label{last2}(\gg^\ell\dot\ot \aa_\ell)\op(\ss^\ell\dot\ot
\bb_\ell)\op(\v^\ell\dot\ot\cc_\ell)\op \dd_\ell\end{equation}
where  $\dd_\ell:=\varphi(E).$
 In what follows using
\cite[Thm. 2.48]{ABG2}, for $\mu=\ell,n,$ we  give the algebraical structure of
$\LL^{^\mu}$ in terms of the ingredients involved in the decomposition of $\LL^\mu$ into finite dimensional irreducible $\gg^\mu-$modules.  Set  $\fa_\mu:=\aa_\mu\op\bb_\mu.$
Then there are a bilinear map
$\cdot_\mu:\fa_\mu\times\fa_\mu\longrightarrow\fa_\mu$
 and a linear map
$*_\mu:\fa_\mu\longrightarrow\fa_\mu$ such that $(\fa_\mu,\cdot_\mu)$
is a unital  associative algebra and $*_\mu$ is an involution on $\fa_\mu$ with
$*_\mu$-fixed points $\aa_\mu$ and $*_\mu$-skew fixed points
$\bb_\mu.$  Also there  is a bilinear map $\cdot_\mu:\fa_\mu\times
\cc_\mu\longrightarrow\cc_\mu$ such that $(\cc_\mu,\cdot_\mu)$ is
a left unital associative $\fa_\mu$-module equipped with a  skew-hermitian form
$f_\mu:\cc_\mu\times\cc_\mu\longrightarrow \fa_\mu.$ Take
$\fb_\mu:=\fb(\fa_\mu,*_\mu,\cc_\mu,f_\mu)$ to be defined as in
Subsection \ref{subsect2-1} and set  $\cdot_\mu,$ $\circ_\mu,$
$[\cdot,\cdot]_\mu,$ $\diamond_\mu$ and $\heart_\mu$ to be the
corresponding features as $\cdot,$ $\circ,$ $[\cdot,\cdot],$
$\diamond$ and $\heart$ defined in Subsection \ref{subsect2-1}.
Also for $\b,\b'\in\fb_\mu,$ set $d_{\b,\b'}^\mu:=d_{\b,\b'}^{\mu,\fb_\mu}.$
 By  \cite[Theorems 2.48, 5.34]{ABG2}, $\dd_\mu$ is a subalgebra of $\LL^\mu$ and there is a subspace $\kk_\mu$ of   the full skew-dihedral homology group $$FH(\fb_\mu)=\{\sum_{i}\{\b_i,\b'_i\}_\mu\mid \sum_{i}d^\mu_{\b_i,\b'_i}=0\}$$ of $\fb_\mu$ such that  $\dd_\mu$ is isomorphic to the quotient algebra  $\{\fb_\mu,\fb_\mu\}_\mu/\kk_\mu.$
For $\b_1,\b_2,$ take $\la\b_1,\b_2\ra_\mu$ to be the element of $\dd_\mu$ corresponding to  $\{\b_1,\b_2\}_\mu+\kk_\mu$, then  one has
$\la\aa_\mu,\bb_\mu\ra_\mu=\la\aa_\mu,\cc_\mu\ra_\mu=\la\bb_\mu,\cc_\mu\ra_\mu=\{0\}$
and
$\dd_\mu=\la\aa_\mu,\aa_\mu\ra_\mu+\la\bb_\mu,\bb_\mu\ra_\mu+\la\cc_\mu,\cc_\mu\ra_\mu.$
Moreover the Lie bracket  on $\LL^\mu$ which is an extension
of the Lie bracket on $\dd_\mu$ is given by
\begin{equation}\label{probc-gen-mu}
\begin{array}{l}
\;[x\ot a,y\ot a']=[x,y]\ot\frac{1}{2}(a\circ_\mu a')+ (x\circ_\mu y)\ot\frac{1}{2}[a,a']_\mu+tr(xy)\la a,a'\ra_\mu,\vspace{1mm}\\
\;[x\ot a,s\ot b]= (x\circ_\mu s)\ot\frac{1}{2}[a,b]_\mu+[x,s]\ot\frac{1}{2}(a\circ_\mu b)=-[s\ot b,x\ot a],\vspace{1mm}\\
\;[s\ot b,t\ot b']=[s,t]\ot\frac{1}{2}(b\circ_\mu b')+ (s\circ_\mu t)\ot\frac{1}{2}[b,b']_\mu+tr(st)\la b,b'\ra_\mu,\vspace{1mm}\\
\;[x\ot a,u\ot c]=xu\ot a\cdot_\mu c=-[u\ot c,x\ot a],\vspace{1mm}\\
\;[s\ot b,u\ot c]=su\ot b\cdot_\mu c=-[u\ot c,s\ot b],\vspace{1mm}\\
\;[u\ot c,v\ot c']=(u\circ v)\ot (c\diamond_\mu c')+ [u, v]_\mu\ot (c\heart_\mu c')+(u,v)\la c,c'\ra_\mu,\vspace{1mm}\\
\;[\la\b,\b'\ra_\mu,x\ot a]=x\ot d^\mu_{\b,\b'}(a)=-[x\ot a,\la\b,\b'\ra_\mu],\vspace{1mm}\\
\;[\la\b,\b'\ra_\mu,s\ot b]=s\ot d^\mu_{\b,\b'}(b)=-[s\ot b,\la\b,\b'\ra_\mu],\vspace{1mm}\\
\;[\la\b,\b'\ra_\mu,u\ot c]=u\ot d^\mu_{\b,\b'}(c)=-[u\ot
c,\la\b,\b'\ra_\mu],\vspace{1mm}\\
\;[\la\b_1,\b_2\ra_\mu,\la\b'_1,\b'_2\ra_\mu]=\la d^\mu_{\b_1,\b_2}(\b'_1),\b'_2\ra_\mu+\la\b'_1,d^\mu_{\b_1,\b_2}(\b'_2)\ra_\mu,
\end{array}
 \end{equation}
for $x,y\in\gg^\mu,$  $s,t\in\ss^\mu,$ $u,v\in\v^\mu,$
$a,a'\in\aa_\mu,$ $b,b'\in\bb_\mu,$ $c,c'\in\cc_\mu$ and
$\b,\b'\in\fb_\mu.$

\begin{Lemma}\label{divide2} We have $\i=\i_0,$ $\jj=\jj_0$ and $\T=\T_0.$\end{Lemma}
\pf It follows from (\ref{probc-gen-mu}), (\ref{last1}) and (\ref{last2}) that
$$
\LL_\a=\left\{
\begin{array}{ll}
\v^n_\a\dot\ot \cc_n &\hbox{if $\a\in R_{sh}$} \\
(\gg^n_\a\dot\ot\aa_n )\op(\ss^n_\a\dot\ot\bb_n )& \hbox{if $\a\in R_{lg}$}\\
\gg^n_{\a}\dot\ot\aa_n &\hbox{if $\a\in R_{ex}$}
\end{array}
\right.$$
and$$
(\LL^{^S})_\a=\left\{
\begin{array}{ll}
(\v^\ell)_\a\dot\ot\cc_\ell&\hbox{if $\a\in S_{sh}$} \\
(\gg^\ell_\a\dot\ot\aa_\ell)\op(\ss^\ell_\a\dot\ot\bb_\ell)& \hbox{if $\a\in S_{lg}$}\\
\gg^\ell_{\a}\dot\ot\aa_\ell&\hbox{if $\a\in S_{ex}.$}
\end{array}
\right.
$$

Now fix  $\a\in S_{ex},$ then
$$(\gg^\ell)_{\a}\dot\ot\aa_\ell=(\LL^{^S})_\a=\LL_\a=\gg^n_{\a}\dot\ot\aa_n .$$
This together with the fact that  $\gg^\ell_{\a}=\gg^n_{\a}$ is a one dimensional
vector space, implies that the vector space $\aa_\ell$ equals the vector space  $\aa_n. $ In particular we  get  $\i=\i_0.$
Next  fix $\a\in S_{sh},$ then we have
$$\v^\ell_{\a}\dot\ot\cc_\ell=\LL^{^S}_\a=\LL_\a=\v^n_{\a}\dot\ot\cc_n. $$ This as above,  implies that
$\T=\T_0.$ Finally  fix $\a\in S_{lg},$ then
$$(\gg^\ell_{\a}\dot\ot\aa_\ell)\op(\ss^\ell_\a\dot\ot\bb_\ell)=\LL^{^S}_\a=\LL_\a=
(\gg^n_{\a}\dot\ot\aa_n )\op(\ss^n_\a\dot\ot\bb_n ).$$
Now as $\ss^\ell_\a=\ss^n_\a$ is a one dimensional
vector space, $\gg^\ell_{\a}=\gg^n_{\a},$ $\bb_\ell\sub\bb_n$ and  $\aa_\ell=\aa_n ,$  we get  that the two vector spaces $\bb_\ell$ and $\bb_n $ are equal and so   $\jj=\jj_0.$\qed

\medskip

As we have already seen,  on the  vector space level, we have   $$\aa_\ell=\aa_n, \bb_\ell=\bb_n,
\cc_\ell=\cc_n$$ which in turn implies that the vector space $\fb_\ell$ equals the vector space  $\fb_n.$ In the following lemma, we show, in addition,  that the algebras $\fb_\ell$ and $\fb_n$ have the same algebraic structure.

\begin{Lemma}\label{divide3} The algebraic structure of $\fb_n$ coincides with the algebraic structure of  $\fb_\ell.$
\end{Lemma}
\pf Using Lemma \ref{divide2}, we set
\begin{equation}
\label{end8}
\aa:=\aa_\ell=\aa_n,\;\;\bb:=\bb_\ell=\bb_n,\;\;\cc:=\cc_\ell=\cc_n .
\end{equation}

Suppose that $i,j,k$ are distinct elements of $I_0.$ Take
$x:=e_{i,j}-e_{\bar j,\bar i}\in\gg^\ell$  and $y:=e_{ j,k}-e_{\bar k,\bar
j}\in\gg^\ell,$ then $$tr(xy)=0,\;[x,y]=e_{i,k}-e_{\bar k,\bar i},\;x\circ_n  y=x\circ_\ell y=e_{i,k}+e_{\bar k,\bar i}.$$ Now  for
$a,a'\in\aa ,$ by (\ref{probc-gen-mu}), we have
{\small\begin{eqnarray*}
[x,y]\ot\frac{1}{2}(a\circ_n  a')+ (x\circ_n y)\ot\frac{1}{2}[a,a']_n
\hspace{-.2cm}&=&\hspace{-.2cm}[x\ot a,y\ot
a']\\
\hspace{-.2cm}&=&\hspace{-.2cm}[x,y]\ot\frac{1}{2}(a\circ_\ell a')+ (x\circ_\ell
y)\ot\frac{1}{2}[a,a']_\ell.
\end{eqnarray*}}
This in turn implies that $$[x,y]\ot(\frac{1}{2}(a\circ_n
a')-\frac{1}{2}(a\circ_\ell a'))= (x\circ_n
y)\ot(\frac{1}{2}[a,a']_\ell-\frac{1}{2}[a,a']_n),$$ but the left hand
side is an element of $\gg\ot \aa$ and the right hand side is an
element of $\ss\ot \bb.$ Therefore as $[x,y]\neq 0$ and $x\circ_n
y\neq 0,$ we get that
$$\frac{1}{2}[a,a']_\ell-\frac{1}{2}[a,a']_n =0\andd \frac{1}{2}(a\circ_n
a')-\frac{1}{2}(a\circ_\ell a')=0.$$
This now implies that
\begin{equation}\label{gen-bc1}
a\cdot_\ell  a'=a\cdot_n  a';\;\; a,a'\in \aa.
\end{equation}

Next take $i$ and $ j$ to be two distinct elements of $I_0.$ Set
$s:=e_{i,j}+e_{\bar j,\bar i}\in\ss^\ell$ and $x:=e_{j,\bar
j}\in\gg^\ell,$ then we have  $$tr(xs)=0,\;[x,s]=e_{j,\bar
i}-e_{i,\bar j},\;x\circ_\ell s=x\circ_n  s=e_{j,\bar i}+e_{i,\bar
j}.$$
 Now for $a\in\aa$ and $b\in\bb,$ by (\ref{probc-gen-mu}), we have
{\small\begin{eqnarray*}
(x\circ_\ell s)\ot \frac{1}{2}[a,b]_\ell+[x,s]\ot \frac{1}{2}
(a\circ_\ell
b)&=&[x\ot a,s\ot b]\\
&=&(x\circ_n  s)\ot \frac{1}{2}[a,b]_n +[x,s]\ot \frac{1}{2}
(a\circ_n  b).
\end{eqnarray*}}
This implies that $$(x\circ_\ell s)\ot(
\frac{1}{2}[a,b]_\ell-\frac{1}{2}[a,b]_n )=[x,s]\ot (\frac{1}{2}
(a\circ_n  b)-\frac{1}{2} (a\circ_\ell b)).$$ Now as before,
one gets that
$$\frac{1}{2}[a,b]_\ell-\frac{1}{2}[a,b]_n =0\andd \frac{1}{2}
(a\circ_n  b)-\frac{1}{2} (a\circ_\ell b)=0.$$
This in particular implies that
\begin{equation}\label{coin-bc2}
a\cdot_\ell b=a\cdot_n  b\andd b\cdot_\ell  a=b\cdot_n  a;\;\;\;\;\;\;\;\;(a\in\aa,\;
b\in\bb).
\end{equation}

Finally, for  distinct fixed elements $i,j,k$ of $I_0,$ set
$s:=e_{i,\bar j}-e_{j,\bar i},t=e_{\bar i,k}-e_{\bar k,i}\in\ss^\ell.$
Then $$tr(st)=0,\;[s,t]=-e_{j,k}+e_{\bar k,\bar j},\;s\circ_\ell
t=s\circ_n  t=-e_{j,k}-e_{\bar k,\bar j}.$$ Therefore for
$b,b'\in\bb ,$ by (\ref{probc-gen-mu}), we have
{\small\begin{eqnarray*}
([s,t]\ot\frac{1}{2}(b\circ_\ell b'))+ ((s\circ_\ell t)\ot
\frac{1}{2}[b,b']_\ell)\hspace{-2mm}&=&\hspace{-2mm}[s\ot b,t\ot
b']\\\hspace{-2mm}&=&\hspace{-2mm}([s,t]\ot\frac{1}{2}(b\circ_n  b'))+ ((s\circ_n  t)\ot
\frac{1}{2}[b,b']_n ).
\end{eqnarray*}}
This implies that $$[s,t]\ot(\frac{1}{2}(b\circ_\ell
b')-\frac{1}{2}(b\circ_n  b'))=(s\circ_n  t)\ot (\frac{1}{2}[b,b']_n
-\frac{1}{2}[b,b']_\ell).$$
Therefore we get $$\frac{1}{2}(b\circ_\ell
b')-\frac{1}{2}(b\circ_n  b')=0\andd \frac{1}{2}[b,b']_n
-\frac{1}{2}[b,b']_\ell=0$$ and so we have
$b\cdot_\ell b'=b\cdot_n  b'$ which together with
(\ref{gen-bc1}) and (\ref{coin-bc2}) implies that
\begin{equation}\label{gen-bc3}
\parbox{4in}{\begin{center}$\fa_\ell=\fa_n $ (as two
algebras).\end{center}}
\end{equation}

Now take $x\in\gg^\ell,$ $s\in\ss^\ell$  and
$u,v\in\v^\ell$ to be such that $xu\neq0$ and $sv\neq0.$ Then for $a\in\aa,$  $b\in\bb$ and $c\in\cc,$ we get using (\ref{probc-gen-mu}) that
\begin{eqnarray*}
xu\ot a\cdot_\ell c&=&[x\ot a, u\ot c]=xu\ot a\cdot_n  c,\\
sv\ot b\cdot_\ell c&=&[s\ot b,v\ot c]=sv\ot b\cdot_n  c.
\end{eqnarray*}
This  implies that   $$a\cdot_\ell  c=a\cdot_n  c\andd b\cdot_\ell  c=b\cdot_n  c$$ for $a\in\aa,$  $b\in\bb$ and $c\in\cc.$
Therefore we have
\begin{equation}\label{gen-coin7}
\cc_\ell=\cc_n  \hbox{ (as two $\fa_n-$modules). }
\end{equation}
Now  we are done thanks to (\ref{gen-bc3}) and (\ref{gen-coin7}).
\qed

\bigskip

Now using Lemma \ref{divide3}, we set $$\fa:=\fa_n=\fa_\ell\andd \fb:=\fb_\ell=\fb_n.$$ Also for $\b,\b'\in\fb,$ we take \begin{equation}\label{setare}\begin{array}{l}\;\b\cdot\b':=\b\cdot_n\b'=\b\cdot_\ell\b',\\
\;[\b,\b']:=\b\cdot\b'-\b'\cdot\b,\\
\;\b\circ\b':=\b\cdot\b'+\b'\cdot\b.\end{array}\end{equation}

\begin{Lemma}\label{divide4}
For $a,a'\in\aa,$ and  $b,b'\in\bb,$ we have $$\begin{array}{l}\la
a,a'\ra_n =(\frac{-1}{\ell}Id_{_{\v^{^\ell}}}+\frac{1}{n }Id_{_\v}
)\ot\frac{1}{2}[a,a'])+\la a,a'\ra_\ell,\\
\la
b,b'\ra_n =(\frac{-1}{\ell}Id_{_{\v^{^\ell}}}+\frac{1}{n }Id_{_\v}
)\ot\frac{1}{2}[b,b'])+\la b,b'\ra_\ell.
\end{array}$$

Also for $c,c'\in \cc,$ $f_\ell(c,c')=f_n(c,c'),$ $c\diamond_\ell c'=c\diamond_n c'$
and $c\heart_\ell c'=c\heart_n c'.$ Moreover we have
\begin{eqnarray*}\la
c,c'\ra_n &=&((\frac{1}{\ell}Id_{_{\v^{^\ell}}}-\frac{1}{n }Id_{_\v}
)\ot\frac{1}{2} c\heart_n c')+\la c,c'\ra_\ell\\&=&((\frac{1}{\ell}Id_{_{\v^{^\ell}}}-\frac{1}{n }Id_{_\v}
)\ot\frac{1}{2} c\heart_\ell c')+\la c,c'\ra_\ell.\end{eqnarray*}\end{Lemma}
\pf
Fix  $x,y\in\gg^{\ell}$ such that $tr(xy)\neq 0.$ For
$a,a'\in\aa,$ consider (\ref{setare}) and use (\ref{probc-gen-mu}) to get
\begin{eqnarray*}
&&([x,y]\ot\frac{1}{2}(a\circ a'))+((x\circ_n
y)\ot\frac{1}{2}[a,a'])+tr(xy)\la a,a'\ra_n =[x\ot a,y\ot a']=\\
&&([x,y]\ot\frac{1}{2}(a\circ a'))+((x\circ_\ell
y)\ot\frac{1}{2}[a,a'])+tr(xy)\la a,a'\ra_\ell.
\end{eqnarray*}
This implies that
$$(\frac{-tr(xy)}{n }Id_{_\v} \ot\frac{1}{2}[a,a'])+tr(xy)\la
a,a'\ra_n
=(\frac{-tr(xy)}{\ell}Id_{_{\v^{^\ell}}}\ot\frac{1}{2}[a,a'])+tr(xy)\la
a,a'\ra_\ell.$$ Therefore we have
\begin{equation}\label{coincide1}\la
a,a'\ra_n =(\frac{-1}{\ell}Id_{_{\v^{^\ell}}}+\frac{1}{n }Id_{_\v}
)\ot\frac{1}{2}[a,a'])+\la a,a'\ra_\ell;\;\;\;\;\;\; (a,a'\in\aa).\end{equation}

Next fix $s,t\in\ss^\ell$ such that $tr(st)\neq 0,$ then for
$b,b'\in\bb,$ by (\ref{probc-gen-mu}), we have
\begin{eqnarray*}
&&([s,t]\ot \frac{1}{2} b\circ b')+((s\circ_n  t)\ot
\frac{1}{2}[b,b'])+tr(st)\la b,b'\ra_n =[s\ot b,t\ot b']=\\
&&([s,t]\ot \frac{1}{2} b\circ b')+((s\circ_\ell t)\ot
\frac{1}{2}[b,b'])+tr(st)\la b,b'\ra_\ell.
\end{eqnarray*}
This implies that $$(\frac{-tr(st)}{n } Id_{_{\v}} \ot
\frac{1}{2}[b,b'])+tr(st)\la b,b'\ra_n
=(\frac{-tr(st)}{\ell} Id_{_{\v^\ell}}\ot
\frac{1}{2}[b,b'])+tr(st)\la b,b'\ra_\ell$$ which in turn implies
that
\begin{equation}\label{coin-gen1} \la b,b'\ra_n =((\frac{-1}{\ell} Id_{_{\v^\ell}}+\frac{1}{n }Id_{_\v} )\ot
\frac{1}{2}[b,b'])+\la b,b'\ra_\ell;\;\;\;\;\;(b,b'\in\bb) .\end{equation}

Now suppose that $i$ and $j$ are two distinct elements of $I_0.$ Take
$u:=v_i$ and $v:=v_{\bar j},$ then $(u,v)=0$ and so $[u,v]_n=[u,v]
.$ Therefore  for all $c,c'\in\cc,$ by (\ref{probc-gen-mu}),  we have
\begin{eqnarray*}
(u\circ v )\ot (c\diamond_n c')+[u,v]\ot (c\heart_n
c')&=&[u\ot c,v\ot c']\nonumber\\
&=&(u\circ v)\ot (c\diamond_\ell  c')+[u,v] \ot c\heart_\ell
c'.\nonumber
\end{eqnarray*}
But  $u\circ v\in\gg,$ $[u,v] $ is a nonzero element of
$\ss,$  $c\diamond_\ell c',c\diamond_n  c'\in\aa$ and
$c\heart_\ell c',c\heart_n  c'\in \bb,$ so  we get that
\begin{equation}\label{last5}
c\diamond c':=c\diamond_\ell c'=c\diamond_n  c' \andd c\heart
c':=c\heart_\ell c'=c\heart_n  c';\;\;\;\;\;\; (c,c'\in\cc).
\end{equation}
This in turn implies that
\begin{equation}
\label{gen-coin5} f(c,c'):= f_\ell(c,c')=f_n(c,c');\;\;\;\;\; (c,c'\in\cc) .
\end{equation}
Next for an element $i$ of $I,$ take  $u:=v_i$ and $v=v_{\bar i}.$
Then for $c,c'\in\cc,$ by (\ref{probc-gen-mu}), we have
\begin{eqnarray*}
&&(u\circ v\ot c\diamond c')+([u,v]_n\ot c\heart
c')+\la c,c'\ra_n=[u\ot c,v\ot c']=\nonumber\\
&&(u\circ v\ot c\diamond c')+([u,v]_\ell \ot c\heart c')+\la
c,c'\ra_\ell ,\nonumber
\end{eqnarray*}
using which, one concludes that  \begin{equation}\label{gen-coin6} \la c,c'\ra_n =\la
c,c'\ra_\ell+(\frac{1}{2\ell}Id_{_\v}-\frac{1}{2n }Id_{_{\v^\ell}}
)\ot c\heart c';\;\; \;\;\;\;(c,c'\in\cc).
\end{equation}
This completes the proof.\qed
\bigskip

\begin{cor}\label{cor1}
Let $\ell<n$ and suppose that  $t\in\bbbn,$ $a_i,a'_i\in\aa,$ $b_i,b'_i\in\bb$ and $c_i,c'_i\in\cc$ for $1\leq i\leq t.$ Then $\sum_{i=1}^t(\la a_i,a'_i\ra_\ell+\la b_i,b'_i\ra_\ell+\la c_i,c'_i\ra_\ell)=0$ if and only if $\sum_{i=1}^t([a_i,a'_i]+[b_i,b_i']-c_i\heart c_i')=0$  and   $\sum_{i=1}^t(\la a_i,a'_i\ra_n+\la b_i,b'_i\ra_n+\la c_i,c'_i\ra_n)=0.$
\end{cor}
\proof   By Lemma \ref{divide4}, we have
\begin{eqnarray*}
&&\sum_{i=1}^t(\la a_i,a'_i\ra_\ell+\la b_i,b'_i\ra_\ell+\la c_i,c'_i\ra_\ell)=\\
&&\sum_{i=1}^t(\la a_i,a'_i\ra_n+\la b_i,b'_i\ra_n+\la c_i,c'_i\ra_n)-\\
&&(\frac{-1}{\ell}Id_{_{\v^{^\ell}}}+\frac{1}{n }Id_{_\v}
)\ot\frac{1}{2}\sum_{i=1}^t([ a_i,a'_i]+[b_i,b'_i]-c_i\heart c'_i).
\end{eqnarray*}
Now as  $\sum_{i=1}^t(\la a_i,a'_i\ra_n+\la b_i,b'_i\ra_n+\la c_i,c'_i\ra_n)\in\dd_n$ and  $(\frac{-1}{\ell}Id_{_{\v^{^\ell}}}+\frac{1}{n }Id_{_\v}
)\ot\frac{1}{2}\sum_{i=1}^t([ a_i,a'_i]+[b_i,b'_i]-c_i\heart c'_i)\in\ss\ot\bb,$ we are done.\qed

\begin{rem}\label{rem1}
{\rm Consider the decomposition (\ref{last3-fin}) for $\LL^{^S}=\LL^\ell$ into finite dimensional  irreducible $\fg^{^S}-$submodules and the  decomposition of $\LL=\LL^n$ into finite dimensional  irreducible $\fg-$submodules as in Lemma \ref{divide1}, then  contemplating the identification (\ref{identification}), we have using Lemma \ref{divide4} that
$$\LL=(\bigoplus_{i\in\i}\gg_i\op\bigoplus_{j\in\jj}\ss_j\op\bigoplus_{t\in\T}\v_t)+ E.$$ Moreover setting $\la\b,\b'\ra^n:=\varphi^{-1}(\la\b,\b'\ra_n)$ and $\la\b,\b'\ra^\ell:=\varphi^{-1}(\la\b,\b'\ra_\ell)$ for $\b,\b'\in \fb,$ we get that $\{\la\b,\b'\ra^n\mid \b,\b'\in\fb\}$ spans $\dd_n$ and that  $\{\la\b,\b'\ra^\ell\mid \b,\b'\in\fb\}$ spans $E.$
Furthermore, thanks to Corollary \ref{cor1},  for $t\in\bbbn,$ $a_i,a'_i\in\aa,$ $b_i,b'_i\in\bb$ and $c_i,c'_i\in\cc$  ($1\leq i\leq t$), $\sum_{i=1}^t(\la a_i,a'_i\ra^\ell+\la b_i,b'_i\ra^\ell+\la c_i,c'_i\ra^\ell)=0$ if and only if $\sum_{i=1}^t([a_i,a'_i]+[b_i,b_i']-c_i\heart c_i')=0$  and   $\sum_{i=1}^t(\la a_i,a'_i\ra^n+\la b_i,b'_i\ra^n+\la c_i,c'_i\ra^n)=0.$  }
\end{rem}

\begin{Proposition}
\label{divide5} For $e,f\in\gg\cup\ss,$ set
$$e\circ f:=ef+fe-\frac{tr(ef)}{\ell}Id_{_{\v^\ell}}.$$ Also
for $\b_1=a_1+b_1+c_1\in\fb$ and $\b_2=a_2+b_2+c_2\in\fb$ with $a_1,a_2\in\aa,$ $b_1,b_2\in\bb$ and $c_1,c_2\in\cc,$ we recall from (\ref{beta*}) that $\b_{_{\b_1,\b_2}}^*:=[a_1,a_2]+[b_1,b_2]-c_1\heart c_2$ and  set $$\la\b_1,\b_2\ra:=\la\b_1,\b_2\ra_\ell,\;\;\b_1^*=c_1,\;\;\b_2^*=c_2$$
and take  $$\dd:=\hbox{span}\{\la a,a'\ra,\la b,b'\ra,\la
c,c'\ra\mid a,a'\in\aa,\;b,b'\in\bb,c,c'\in\cc\},$$ then contemplating  (\ref{last5}), we  have
\begin{eqnarray*}\LL&=& (\gg\ot \aa)\op(\ss\ot
\bb)\op(\v\ot\cc)\op\la\fb,\fb\ra_n\nonumber\\&=& ((\gg\ot
\aa)\op(\ss\ot \bb)\op(\v\ot\cc))+\dd\end{eqnarray*} with the Lie bracket given by{\small
\begin{equation}\label{probc-fin}
\begin{array}{l}
\;[x\ot a,y\ot a']=[x,y]\ot\frac{1}{2}(a\circ a')+ (x\circ y)\ot\frac{1}{2}[a,a']+tr(xy)\la a,a'\ra,\vspace{1mm}\\
\;[x\ot a,s\ot b]= (x\circ s)\ot\frac{1}{2}[a,b]+[x,s]\ot\frac{1}{2}(a\circ b)=-[s\ot b,x\ot a],\vspace{1mm}\\
\;[s\ot b,t\ot b']=[s,t]\ot\frac{1}{2}(b\circ b')+ (s\circ t)\ot\frac{1}{2}[b,b']+tr(st)\la b,b'\ra,\vspace{1mm}\\
\;[x\ot a,u\ot c]=xu\ot a\cdot c=-[u\ot c,x\ot a],\vspace{1mm}\\
\;[s\ot b,u\ot c]=su\ot b\cdot c=-[u\ot c,s\ot b],\vspace{1mm}\\
\;[u\ot c,v\ot c']=(u\circ v)\ot (c\diamond c')+ [u, v]\ot (c\heart c')+(u,v)\la c,c'\ra,\vspace{1mm}\\
\;[\la \b_1,\b_2\ra,x\ot a]=
\frac{-1}{4\ell}(x\circ
Id_{_{\v^\ell}}\ot[a,\b_{_{\b_1,\b_2}}^*]+[x,Id_{_{\v^\ell}}]\ot a\circ \b_{_{\b_1,\b_2}}^*),\vspace{1mm}\\
\;[\la \b_1,\b_2\ra,\hspace{-1mm}s\ot
b]\hspace{-1mm}=\hspace{-1mm}\frac{-1}{4\ell}([s,Id_{_{\v^\ell}}\hspace{-.5mm}]\hspace{-1mm}\ot\hspace{-1mm} (b\circ \b_{_{\b_1,\b_2}}^*\hspace{-1mm})\hspace{-1mm}+\hspace{-1mm}(s\circ
Id_{_{\v^\ell}})\hspace{-1mm}\ot \hspace{-1mm}[b, \b_{_{\b_1,\b_2}}^*\hspace{-1mm}]\hspace{-1mm}+\hspace{-1mm}2tr(sId_{\v^\ell})\la b,\b_{_{\b_1,\b_2}}^*\hspace{-.5mm}\ra),\vspace{1mm}\\
\;[\la \b_1,\b_2\ra,v\ot
c]=\frac{1}{2\ell}Id_{_{\v^\ell}}v\ot (\b_{_{\b_1,\b_2}}^*\cdot c)-\frac{1}{2}v\ot
(f(c,\b^*_2)\cdot \b^*_1+f(c,\b^*_1)\cdot \b^*_2)\\
\;[\la\b_1,\b_2\ra,\la\b'_1,\b'_2\ra]=\la d^\ell_{\b_1,\b_2}(\b'_1),\b'_2\ra+\la\b'_1,d^\ell_{\b_1,\b_2}(\b'_2)\ra
\end{array}
 \end{equation}}
for $x,y\in\gg,$  $s,t\in\ss,$ $u,v\in\v,$
$a,a'\in\aa,$ $b,b'\in\bb,$ $c,c'\in\cc,$ $\b_1,\b_2,\b_1',\b'_2\in\fb.$
\end{Proposition}

\pf Suppose that $x,y\in\gg,$  $s,t\in\ss, $  $u,v\in\v, $ $a,a'\in\aa,$  $b,b'\in\bb $ and $c,c'\in\cc,$ then (\ref{probc-gen-mu}) (for $\mu=n$)  together with Lemma \ref{divide4} implies that
{\small\begin{eqnarray}
[x\ot a,y\ot a']&=&([x,y]\ot \frac{1}{2}(a\circ a'))+((x\circ_n
y)\ot \frac{1}{2}[a,a'])+tr(xy)\la a,a'\ra_{n}\nonumber\\
&=&([x,y]\ot \frac{1}{2}(a\circ a'))+((x\circ_n  y)\ot
\frac{1}{2}[a,a'])\nonumber\\&+&tr(xy)(\frac{-1}{\ell} Id_{_{\v^{^\ell}}}+\frac{1}{n
}Id_{_\v} )\ot\frac{1}{2}[a,a'])+tr(xy)\la
a,a'\ra_\ell\label{pro-bc-gen-a}\\
&=&([x,y]\ot \frac{1}{2}(a\circ a'))\nonumber\\&+&(((x\circ_n
y)+tr(xy)(\frac{-1}{\ell}Id_{_{\v^{^\ell}}}+\frac{1}{n } Id_{_\v}
))\ot
\frac{1}{2}[a,a'])+tr(xy)\la a,a'\ra_\ell\nonumber\\
&=&([x,y]\ot \frac{1}{2}(a\circ a'))+((x\circ y)\ot
\frac{1}{2}[a,a'])+tr(xy)\la a,a'\ra.\nonumber
\end{eqnarray}}
Similarly we have
$$[s\ot b,t\ot b']=([s,t]\ot \frac{1}{2}(b\circ b'))+((s\circ t)\ot
\frac{1}{2}[b,b'])+tr(st)\la b,b'\ra$$ and $$[u\ot c,v\ot
c']=((u\circ v)\ot c\diamond c')+([u,v]\ot (c\heart c'))+(u,v)\la
c,c'\ra.$$
%

Now for $a_1,a_2\in\aa,$ $b_1,b_2\in\bb$ and $c_1,c_2\in\cc,$ set
$b^*:=[a_1,a_2]+[b_1,b_2]-c_1\heart c_2]$ and  take
$s_n:=(1/\ell)Id_{_{\v^\ell}}-(1/n)Id_{_{\v}}.$ Then for  $x\in\gg,$ $s\in\ss,$ $v\in\v,$
 $a\in\aa,$ $b\in \bb,$ and $c\in\cc,$
one can see that {\small $$[x\ot a,s_n\ot
b^*]=\frac{1}{2\ell}(x\circ
Id_{_{\v^\ell}}\ot[a,b^*]+[x,Id_{_{\v^\ell}}]\ot a\circ
b^*)-\frac{1}{n}x\ot[a,b^*],$$}

{\small $$[s\ot b,s_n\ot
b^*]=\frac{1}{2\ell}([s,Id_{_{\v^\ell}}]\ot b\circ b^*+(s\circ
Id_{_{\v^\ell}})\ot [b, b^*])-\frac{1}{n}s \ot
[b, b^*]\\
+tr(ss_n)\la b,b^*\ra$$ and
$$[s_n\ot b^*,v\ot c]=
\frac{1}{\ell}Id_{_{\v^\ell}}v\ot b^*\cdot c-\frac{1}{n}v\ot
b^*\cdot c.$$}
%
%
%

We next note that {\small \begin{eqnarray*}[\la a_1,a_2\ra_n+\la
b_1,b_2\ra_n+\la c_1,c_2\ra_n,x\ot a]&=&x\ot
(d^n_{a_1,a_2}+d^n_{b_1,b_2}+d^n_{c_1,c_2})(a)\\
&=&\frac{1}{2n}x\ot[b^*,a]\\
\;[\la a_1,a_2\ra_n+\la b_1,b_2\ra_n+\la c_1,c_2\ra_n,s\ot
b]&=&s\ot
(d^n_{a_1,a_2}+d^n_{b_1,b_2}+d^n_{c_1,c_2})(b)\\
&=&\frac{1}{2n}s\ot[b^*,b]\\
\;[\la a_1,a_2\ra_n+\la b_1,b_2\ra_n+\la c_1,c_2\ra_n,v\ot
c]&=&v\ot
(d^n_{a_1,a_2}+d^n_{b_1,b_2}+d^n_{c_1,c_2})(c)\\
&=&\frac{1}{2n}v\ot b^*\cdot c\\&-&v\ot \frac{1}{2}(f(c,c_2)\cdot
c_1+f(c,c_1)\cdot c_2).
\end{eqnarray*}}
Therefore using Lemma \ref{divide4}, an easy verification gives
that {\small $$[\la a_1,a_2\ra+\la b_1,b_2\ra+\la c_1,c_2\ra,x\ot
a]=-\frac{1}{4\ell}(x\circ
Id_{_{\v^\ell}}\ot[a,b^*]+[x,Id_{_{\v^\ell}}]\ot a\circ b^*),$$

\begin{eqnarray*}[\la
a_1,a_2\ra+\la b_1,b_2\ra+\la c_1,c_2\ra,s\ot
b]&=&-\frac{1}{4\ell}([s,Id_{_{\v^\ell}}]\ot b\circ b^*+(s\circ
Id_{_{\v^\ell}})\ot [b, b^*])\\&-&\frac{1}{2\ell}tr(sId_{\v^\ell})\la b,b^*\ra,\end{eqnarray*} and $$[\la a_1,a_2\ra+\la
b_1,b_2\ra+\la c_1,c_2\ra,v\ot
c]=\frac{1}{2\ell}Id_{_{\v^\ell}}v\ot b^*\cdot c-\frac{1}{2}v\ot
(f(c,c_2)\cdot c_1+f(c,c_1)\cdot c_2).$$}
These   together with (\ref{probc-gen-mu})
complete the proof.\qed

\subsubsection{\textbf{Types $A$ and $D$}}  Suppose that $I$ is an index set of cardinality $n+1>5$ and  $I_0$ is a subset of $I$ of cardinality $\ell+1>5.$ Suppose that $R$ is an irreducible finite root system  of type $X=\dot A_I$ or $D_I.$ Suppose that $\v$ is a vector space with a basis $\{v_i\mid i\in I\}$ and take $\gg$ to be the  finite dimensional split simple Lie algebra of type $X$ as in   Lemma \ref{type-a-alg} or Lemma \ref{type-d-alg} respectively and set $\gg^\ell:=\gg_{_{I_0}}.$ Suppose that  $\v^\ell$ is the subspace of $\v$ spanned by $\{v_i\mid i\in I_0\}.$ We  take $Id_{_\v}$ to be the identity map on $\v$ and define $Id_{_{\v^\ell}}$ as follows:
$$\begin{array}{c}Id_{_{\v^{\ell}}}:\v\longrightarrow \v\\
v_i\mapsto v_i,\;v_{j}\mapsto 0;\;\;
(i\in I_0,\; j\in I\setminus I_0).
\end{array}$$

\begin{Theorem}\label{type-a}
Suppose that $\LL$ is a Lie algebra graded by the irreducible finite root system $R$ of type $X=\dot A_I$ or $D_I$ with grading pair $(\fg,\fh)$ and let  $S$ be the  irreducible full  subsystem of $R$ of type $\dot A_{I_0}$ or $D_{I_0}$ respectively. 

(i) Consider $\LL^{^S}$ as a $\fg^{^S}-$module and  take
\begin{equation}\label{last3-a}\LL^{^S}=\bigoplus_{i\in
\i}\fg_i\op E\end{equation} to be the decomposition of $\LL^{^S}$
into finite dimensional irreducible $\fg^{^S}-$ submodules  in which
$\i$ is an  index set and for
$i\in\i,$  $\fg_i$ is isomorphic to
$\fg^{^S}(\simeq\gg^\ell),$
 and $E$ is a trivial
$\fg^{^S}-$submodule.
Then there exists a class $\{\dd_n,\gg_i\mid i\in \i\}$ of  finite dimensional   $\fg-$submodules of $\LL$ such that
\begin{itemize}
\item   $\dd_n$ is a trivial $\fg-$module and   $\gg_i$  is  isomorphic to $\fg\simeq\gg$   for $i\in\i,$   \\
\item $\fg_i\sub\gg_i,$  ($i\in\i$),\\ \item
$\LL=\bigoplus_{i\in\i}\gg_i\op\dd_n.$
\end{itemize}

 \medskip

(ii) Take $\aa$ to be a vector  space with basis $\{a_i\mid i\in\i\}$and identify $\LL$ with $(\gg\ot \aa)\op\dd_n,$ say via the natural identification  $$\varphi:\LL\longrightarrow (\gg\ot\aa)\op\dd_n .$$ Transfer the Lie algebraic structure of $\LL$ to $(\gg\ot\aa)\op\dd_n.$ Then $\dd_\ell:=\varphi(E)$ is a subalgebra of $\varphi(\LL^{^S})=(\gg^{^\ell}\dot\ot \aa)\op\dd_\ell$  and $\dd_n$ is a subalgebras of $(\gg\ot\aa)\op\dd_n.$ Moreover the vector space  $\aa$ is equipped with an associative  algebraic structure if $X=\dot A_I$ and with a commutative associative algebraic structure if $X=D_I.$

(iii) There is a  subspace $\kk_1$ of  the full skew-dihedral homology group of $\aa$ with respect to $n$ and a subspace  $\kk_2$ of the full skew-dihedral homology group of  $\aa$ with respect to $\ell$ such that $\dd_n$ and $\dd_\ell$ are isomorphic to the quotient  algebras $\{\aa,\aa\}_n/\kk_1$ and $\{\aa,\aa\}_\ell/\kk_2$
respectively, say via $$\psi_1:\{\aa,\aa\}_n/\kk_1\longrightarrow \dd_n\andd\psi_2:\{\aa,\aa\}_\ell/\kk_2\longrightarrow \dd_\ell.$$

(iv) For $a,a'\in\aa,$ take $$\la a,a'\ra_n:=\psi_1(\{a,a'\}_n+\kk_1)\andd\la a,a'\ra_\ell:=\psi_2(\{a,a'\}_\ell+\kk_2).
$$ Then for $a,a'\in \aa,$ we have $$\la a,a'\ra_n=\la a,a'\ra_\ell+((\frac{1}{n+1} Id_{_{\v}}-\frac{1}{\ell+1} Id_{_{\v^\ell}})\ot (aa'-a'a).$$

(v) For $a,a'\in\aa,$ set $$\la a,a'\ra^n:=\varphi^{-1}(\la a,a'\ra_n)\andd\la a,a'\ra^\ell:=\varphi^{-1}(\la a,a'\ra_\ell).$$ If $\ell < n,$ then for $a_1,a'_1,\ldots,a_t,a'_t\in\aa,$ we have $\sum_{i=1}^t\la a_i,a'_i\ra^\ell=0$ if and only if  $\sum_{i=1}^t\la a_i,a'_i\ra^n=0$ and $\sum_{i=1}^t[a_i,a'_i]=0.$

(vi) For $x,y\in \gg,$ set $x\circ y:= xy +yx - \frac{2tr(xy)}{\ell+1}Id_{\v^\ell}$ and for $a,a'\in\aa,$ set $\la a,a'\ra:=\la a,a'\ra_\ell.$ The Lie bracket on $(\gg\ot\aa)\op \dd_n=(\gg\ot\aa)+ \dd_\ell$ is given by
\begin{equation}\label{probc-fin-a}
\begin{array}{l}
\hbox{\small$[x\ot a,y\ot a']$}=\left\{\begin{array}{ll}\hbox{\small$\;[x,y]\ot\frac{1}{2}(a\circ a')+ (x\circ y)\ot\frac{1}{2}[a,a']+tr(xy)\la a,a'\ra$}& \hbox{\small$X=\dot A_I$},\vspace{1mm}
\\\hbox{\small$[x,y]\ot aa'+tr(xy)\la a,a'\ra$}&\hbox{\small$X=D_I,$}\end{array}\right.
\\
\hbox{\small$[\la
a_1,a_2\ra,x\ot a]=$}\left\{\begin{array}{ll}
\hbox{\small$\frac{-1}{2(\ell+1)}(x\circ
Id_{_{\v^\ell}}\ot[a,[a_1,a_2]]$}&\\\hbox{\small$+[x,Id_{_{\v^\ell}}]\ot a\circ [a_1,a_2]+2tr(Id_{_{\v^\ell}}x)\la a,[a_1,a_2]\ra)$},&\hbox{\small$X=\dot A_I,$}\\
0&\hbox{\small$X=D_I,$}\end{array}\right.\vspace{1mm}\\
\hbox{\small$[\la a_1,a_2\ra,\la a'_1,a'_2\ra]$}=\left\{\begin{array}{ll}\hbox{\small$\la d^{\ell,\aa}_{a_1,a_2}(a'_1),a'_2\ra+\la a'_1,d^{\ell,\aa}_{a_1,a_2}(a'_2)\ra,$}& \hbox{\small$X=\dot A_I,$}\\
0&\hbox{\small$X=D_I,$}\end{array}\right.
\end{array}
 \end{equation}
for $x,y\in\gg,$
$a,a',a_1,a_2,a'_1,a'_2\in\aa.$
\end{Theorem}

\medskip

\subsubsection{\textbf{Types $B$ and $C$}}\label{subsub b-c}
Suppose that $I$ is a nonempty  index set of cardinality $n$ greater than 4 and  $I_0$ is a subset of $I$ of cardinality $\ell>4.$ Take $\gg$ to be either $\mathfrak{o}_B(I)$ or $\mathfrak{sp}(I).$ Suppose that $\v$ is a vector space with a basis  $\{v_0,v_i,v_{\bar i}\mid i\in I\}$ equipped with a nondegenerate  symmetric bilinear form  $\fm$ as in (\ref{form-b-alg}) if $\gg=\mathfrak{o}_B(I)$ and it is a vector space with a basis  $\{v_i,v_{\bar i}\mid i\in I\}$      equipped with a nondegenerate skew-symmetric bilinear form  $\fm$ as in (\ref{form-c}) if $\gg:=\mathfrak{sp}(I).$ Consider (\ref{simple-b-mod}) and (\ref{simple-c})
and set $\v^\ell:=\v_{_{I_0}}.$ Set $$J:=\left\{\begin{array}{ll}I_0\cup\bar{I_0}\cup\{0\}& \hbox{if $\gg=\mathfrak{o}_B(I)$}\\
I_0\cup\bar{I_0}&\hbox{if $\gg=\mathfrak{sp}(I)$}\end{array}\right.$$ and define $Id_{_{\v^\ell}}:\v\longrightarrow \v$  to be the linear transformation defined by $$v_i\mapsto \left\{\begin{array}{ll}v_i& \hbox{if $i\in J$}\\
0& \hbox{if $i\in I\cup\bar{I} \setminus J.$}\end{array}\right.
$$
Finally set $\ss:=\v$ and $\ss^\ell:=\v^\ell$ if $\gg:=\mathfrak{o}_B(I)$  and take  $\ss$ and $\ss^\ell:=\ss_{_{I_0}}$ to be as in (\ref{module-s-c}) and (\ref{simple-c}) respectively  if $\gg=\mathfrak{sp}(I).$

\begin{Theorem} \label{type-fini-b}
Suppose that $\LL$ is a Lie algebra graded by a root system $R$ of type $X=B_I $ or  $C_I$ with grading pair $(\fg,\fh)$ and let  $S$ be the irreducible full  subsystem of $R$ of  type $B_{I_0}$ or $C_{I_0}$ respectively.

(i) Consider $\LL^{^S}$ as a $\fg^{^S}-$module and  take
\begin{equation}\label{last3-b}\LL^{^S}=\bigoplus_{i\in
\i}\fg_i\op\bigoplus_{j\in\jj} \sfs_j\op E\end{equation} to be the decomposition of $\LL^{^S}$
into finite dimensional irreducible $\fg^{^S}-$ submodules in which
$\i,\jj$ are  index sets and for
$i\in\i$ and $j\in \jj,$ $\fg_i$ is isomorphic to
$\fg^{^S}(\simeq\gg^\ell),$  $\sfs_j$ is isomorphic to $\ss^\ell,$
 and $E$ is a trivial
$\fg^{^S}-$submodule.
Then there exists a class $\{\dd_n,\gg_i,\ss_j\mid i\in \i,j\in \jj\}$ of  finite dimensional   $\fg-$submodules of $\LL$ such that
\begin{itemize}
\item   $\dd_n$ is a trivial $\fg-$module,  $\gg_i$  is  isomorphic to $\fg(\simeq\gg)$ and  $\ss_j$
 is  isomorphic to $\ss,$   for $i\in\i,j\in \jj,$   \\
\item $\fg_i\sub\gg_i,$   $\sfs_j\sub\ss_j,$
 ($i\in\i,$ $j\in\jj$),\\ \item
$\LL=\bigoplus_{i\in\i}\gg_i\op\bigoplus_{j\in\jj}\ss_j\op\dd_n.$
\end{itemize}

 \medskip

(ii) Take $\aa$ and $\bb$ to be vector spaces with bases $\{a_i\mid i\in\i\}$ and $\{b_j\mid j\in\jj\}$ respectively and identify $\LL$ with $(\gg\ot \aa)\op(\v\ot \bb)\op\dd_n,$ say via the natural identification  $$\varphi:\LL\longrightarrow \bigoplus_{i\in\i}\gg_i\op\bigoplus_{j\in\jj}\ss_j\op\dd_n .$$ Transfer the Lie algebraic structure of $\LL$ to $\bigoplus_{i\in\i}\gg_i\op\bigoplus_{j\in\jj}\ss_j\op\dd_n.$ Then $\dd_\ell:=\varphi(E)$ and $\dd_n$ are subalgebras of $\LL.$

\medskip

(iii) Set  $\fa:=\aa\op\bb.$ If $\gg=\mathfrak{o}_B(I),$ $\aa$ is equipped with a unital commutative associative algebraic structure and the vector space $\bb$ is equipped with  a unital  $\aa-$module structure. Also there is a symmetric $\aa-$bilinear form $f:\bb\times \bb\longrightarrow \aa$ and  $\fa=\mathcal{J}(f,\bb).$ 
Also if $\gg=\mathfrak{sp}(I),$ $\fa$
is equipped with a star algebraic structure with an involution $*$ such that $\aa$ (resp. $\bb$) is the set of $*-$fixed (resp. $*$-skew fixed) points of $\fa.$

\medskip

(iv) There is a  subspace $\kk_1$ of the full skew-dihedral homology group of $\fa$ with respect to $n$ and  a subspace $\kk_2$ of  the full skew-dihedral homology group of $\fa$  with respect to $\ell$ such that $\dd_n$ and $\dd_\ell$ are isomorphic to the quotient  algebras $\{\fa,\fa\}_n/\kk_1$ and $\{\fa,\fa\}_\ell/\kk_2$
respectively, say via $\psi_1:\{\fa,\fa\}_n/\kk_1\longrightarrow \dd_n$ and $\psi_2:\{\fa,\fa\}_\ell/\kk_2\longrightarrow \dd_\ell.$

(v) For $\a,\a'\in\fa,$ take $$\la \a,\a'\ra_n:=\{\a,\a'\}_n+\kk_1\andd\la \a,\a'\ra_\ell:=\{\a,\a'\}_\ell+\kk_2.
$$ Then if $\a,\a'\in \aa$ or $\a,\a'\in\bb,$ we have  $$\la \a,\a'\ra_n=\la \a,\a'\ra_\ell+((\frac{1}{n} Id_{_\v}-\frac{1}{\ell} Id_{_{\v^\ell}})\ot (1/2)(\a\a'-\a'\a)).$$

(vi) For $\a,\a'\in\fa,$ set $$\la \a,\a'\ra^n:=\varphi^{-1}(\la \a,\a'\ra_n)\andd\la \a,\a'\ra^\ell:=\varphi^{-1}(\la \a,\a'\ra_\ell).$$ If $\ell < n,$ then for $a_1,a'_1,\ldots,a_t,a'_t\in\aa$ and  $b_1,b'_1,\ldots,b_t,b'_t\in\bb,$ we have $\sum_{i=1}^t\la a_i,a'_i\ra^\ell+\sum_{i=1}^t\la b_i,b'_i\ra^\ell=0$ if and only if  $$\sum_{i=1}^t\la a_i,a'_i\ra^n+\sum_{i=1}^t\la b_i,b'_i\ra^n=0 \andd \sum_{i=1}^t[a_i,a'_i]+\sum_{i=1}^t[b_i,b'_i]=0.$$

(vii) For $e,f\in\gg\cup\ss,$ set
$$e\circ f:=ef+fe-\frac{tr(ef)}{\ell}Id_{_{\v^\ell}},$$ and for $\a,\a'\in\fa,$ set $$\la\a,\a'\ra:=\la\a,\a'\ra_\ell,$$  the Lie bracket on $(\gg\ot\aa)\op (\ss\ot \bb) \op \dd_n=((\gg\ot\aa)\op (\ss\ot \bb))+ \dd_\ell$ is given by
\begin{equation}\label{probc-fin-b}
\begin{array}{l}
\;[x\ot a,y\ot a']=[x,y]\ot a a'+ tr(xy)\la a,a'\ra,\vspace{1mm}\\
\;[x\ot a,s\ot b]=xs\ot ab,  \vspace{1mm}\\
\;[s\ot b,t\ot b']=D_{s,t}\ot f(b,b')+(s,t)\la b,b'\ra, \vspace{1mm}\\
\;[\la\a,\a'\ra,x\ot a]=x\ot d^{\ell,\fa}_{\a,\a'}(a),\vspace{1mm}\\
\;[\la\a,\a'\ra,s\ot b]=s\ot d^{\ell,\fa}_{\a,\a'}(b),\vspace{1mm}\\
\;[\la \a_1,\a_2\ra,\la \a'_1,\a'_2\ra]=\la d^{\ell,\fa}_{\a_1,\a_2}(\a'_1),\a'_2\ra+\la\a'_1,d^{\ell,\fa}_{\a_1,\a_2}(\a'_2)\ra.
\end{array}
 \end{equation}
(see Definition \ref{yoshii}) for $x,y\in\gg,$ $s,t\in\ss,$
$a,a'\in\aa,$ $b,b'\in\bb,$ and $\a,\a',\a_1,\a_2,$ $\a'_1,\a'_2\in\fa$ if $\gg=\mathfrak{o}_B(I)$ and it is given by
\begin{equation}\label{probc-fin-c}
\hbox{\small$\begin{array}{l}
\;[x\ot a,y\ot a']=[x,y]\ot\frac{1}{2}(a\circ a')+ (x\circ y)\ot\frac{1}{2}[a,a']+tr(xy)\la a,a'\ra,\vspace{1mm}\\
\;[x\ot a,s\ot b]=(x\circ s)\ot\frac{1}{2}[a,b]+[x,s]\ot \frac{1}{2} (a\circ b),  \vspace{1mm}\\
\;[s\ot b,t\ot b']=[s,t]\ot\frac{1}{2}(b\circ b') + (s\circ t)\ot \frac{1}{2}[b,b']+tr(st)\la b,b'\ra, \vspace{1mm}\\
\;[\la \a,\a'\ra,x\ot a]=
\frac{-1}{4\ell}((x\circ
Id_{_{\v^\ell}})\ot[a,\b_{\a,\a'}^*]+[x,Id_{_{\v^\ell}}]\ot (a\circ \b_{\a,\a'}^*)),\vspace{1mm}\\
\;[\la \a,\a'\ra,s\ot b]\hspace{-1mm}=\hspace{-1mm}\frac{-1}{4\ell}([s,Id_{_{\v^\ell}}]\hspace{-1mm}\ot\hspace{-1mm} (b\circ \b_{\a,\a'}^*)\hspace{-1mm}+\hspace{-1mm}(s\circ
Id_{_{\v^\ell}})\hspace{-1mm}\ot\hspace{-1mm} [b, \b_{\a,\a'}^*]+2tr(sId_{_{\v^\ell}})\la b,\b^*_{\a,\a'}\ra),\vspace{1mm}\\
\;[\la \a_1,\a_2\ra,\la \a'_1,\a'_2\ra]=\la d^{\ell,\fa}_{\a_1,\a_2}(\a'_1),\a'_2\ra+\la\a'_1,d^{\ell,\fa}_{\a_1,\a_2}(\a'_2)\ra.
\end{array}$}
 \end{equation}
(see (\ref{beta*})) for $x,y\in\gg,$ $s,t\in\ss,$
$a,a'\in\aa,$ $b,b'\in\bb,$ $\a,\a',\a_1,\a_2,\a'_1,\a'_2\in\fa$ if $\gg=\mathfrak{sp}(I).$

\end{Theorem}

\section{Root graded Lie algebras - general case}
In this section, we discuss certain  recognition theorems  to characterize    Lie algebras graded by an infinite   irreducible  locally finite root system.  The main target of the present section is  to generalize the  decomposition (\ref{*}) for Lie algebras graded by infinite root systems. For a Lie algebra $\LL$ graded by an infinite locally finite root system with grading pair $(\fg,\fh)$,  we first   decompose $\LL$ as a direct sum of a certain  subalgebra of $\LL$  and a  certain locally finite completely reducible  $\fg-$submodule. This in particular  results in a generalized  decomposition for $\LL$  as in  (\ref{*}). We next reconstruct the structure of $\LL$ in terms of the ingredients  involved in this decomposition.  Moreover we prove that any Lie algebra graded by an irreducible locally finite root system arises in this way. As in the previous section, we concentrate our attention on type $BC$ and for other types, we just report the results.
\subsection{Recognition theorem for  type $BC$}
Suppose that  $I$ is an infinite index set and $\ell $ is a
positive integer greater than 3. We assume     $R$ is   an irreducible
locally finite root system of type $BC_I$ and take  $\gg,$ $\ss$ and $\v$ to be
as in Lemmas \ref{type-c-alg} and \ref{rep-local}. We show that an $R-$graded Lie algebra $\LL$ can be decomposed into \begin{equation}\label{**}(\gg\ot \aa)\op(\ss\ot\bb)\op(\v\ot \cc)\op\dd\end{equation} in which $\aa,\bb$ and $\cc$  are vector spaces and $\dd$ is a subalgebra of $\LL.$ We equip $\fb:=\aa\op\bb\op\cc$ with a unital  associative star algebraic structure and show that $\dd$ can be expressed as a quotient of the algebra $\{\fb,\fb\}_\ell$ by a subspace of  the full skew-dihedral homology group  of $\fb$ with respect to $\ell$ satisfying the uniform  property on $\fb.$ Conversely, for vector spaces $A,B,C$ and $D$ with  specific natures, we  form the decomposition (\ref{**}), equip it with a Lie bracket and show that it is an $R-$graded Lie algebra.
%
%

\begin{Theorem}\label{typebc}
Suppose that  $I$ is an infinite index set and  $\ell$ is an integer greater than 3. Assume  $R$ is an irreducible  locally finite root system of type $BC_I$ and  $\v$ is a vector space with a basis $\{v_i\mid i\in I\cup\bar I\}.$ Suppose that $\fm$ is a bilinear form as in (\ref{form-c}), set $\gg:=\mathfrak{sp}(I)$ and consider $\ss$ as in (\ref{module-s-c}). Fix a subset $I_0$ of $I$ of cardinality $\ell$ and  take $R_0$ to be the  full irreducible  subsystem  of $R$ of type $BC_{I_0}.$ Suppose that  $\{R_\lam\mid\lam\in\Lam\}$ is the class  of all finite irreducible full subsystems of $R$
containing $R_0,$ where $\Lam$ is an index set containing zero.
For $\lam\in\Lam,$ take $\gg^\lam$ as in Lemma \ref{simple-c-alg} and  $\v^\lam,\ss^\lam$ as in (\ref{simple-c}), also define $$\begin{array}{c}\mathfrak{I}_\lam:\v\longrightarrow\v\\
v_i\mapsto\left\{\begin{array}{ll}v_i& i\in I_\lam\cup \bar I_{\lam}\\
0&\hbox{otherwise.}\end{array}\right.\end{array}$$
For $e,f\in\gg\cup\ss,$ define $$e\circ f:=ef+fe-\frac{tr(ef)}{l}\mathfrak{I}_0.$$

(i) Suppose that $(\fa,*,\cc,f)$ is a coordinate quadruple of type $BC$ and $\aa,$ $\bb$ are  $*$-fixed and $*$-skew fixed points of  $\fa$ respectively. Set $\fb:=\fb(\fa,*,\cc,f)$ and take $[\cdot,\cdot], \circ,\heart,\diamond$ to be  as in Subsection \ref{subsect2-1}. For $\b_1,\b_2\in\fb,$ consider $d_{\b_1,\b_2}^{\ell,\fb}$  as in (\ref{derivbc}) and take $\b^*_{\b_1,\b_2},\b_1^*$ and $\b_2^*$ as in Proposition \ref{divide5}. For a subset $\kk$ of $FH(\fb)$ satisfying the uniform property on $\fb,$ set $$\LL(\fb,\kk):=(\gg\ot\aa)\op(\ss\ot \bb)\op(\v\ot\cc)\op(\{\fb,\fb\}_\ell/\kk).$$  Then setting $\la \b,\b'\ra:=\{\b,\b'\}+\kk,$ $\b,\b'\in \fb,$ $\LL(\fb,\kk)$ together with
{\small\begin{equation}\label{probc-gen}
\begin{array}{l}
\;[x\ot a,y\ot a']=[x,y]\ot\frac{1}{2}(a\circ a')+ (x\circ y)\ot\frac{1}{2}[a,a']+tr(xy)\la a,a'\ra,\vspace{1mm}\\
\;[x\ot a,s\ot b]= (x\circ s)\ot\frac{1}{2}[a,b]+[x,s]\ot\frac{1}{2}(a\circ b)=-[s\ot b,x\ot a],\vspace{1mm}\\
\;[s\ot b,t\ot b']=[s,t]\ot\frac{1}{2}(b\circ b')+ (s\circ t)\ot\frac{1}{2}[b,b']+tr(st)\la b,b'\ra,\vspace{1mm}\\
\;[x\ot a,u\ot c]=xu\ot a\cdot c=-[u\ot c,x\ot a],\vspace{1mm}\\
\;[s\ot b,u\ot c]=su\ot b\cdot c=-[u\ot c,s\ot b],\vspace{1mm}\\
\;[u\ot c,v\ot c']=(u\circ v)\ot (c\diamond c')+ [u, v]\ot (c\heart c')+(u,v)\la c,c'\ra,\vspace{1mm}\\
\;[\la \b_1,\b_2\ra,x\ot a]=
\frac{-1}{4\ell}((x\circ
\mathfrak{I}_0)\ot[a,\b_{_{\b_1,\b_2}}^*]+[x,\mathfrak{I}_0]\ot (a\circ \b_{_{\b_1,\b_2}}^*)),\vspace{1mm}\\
\;[\la \b_1,\b_2\ra,s\ot
b]\hspace{-1mm}=\hspace{-1mm}\frac{-1}{4\ell}([s,\mathfrak{I}_0\hspace{-.5mm}]\hspace{-1mm}\ot\hspace{-.5mm}( b\circ \b_{_{\b_1,\b_2}}^*)\hspace{-1mm}+\hspace{-1mm}(s\circ
\mathfrak{I}_0)\hspace{-1mm}\ot \hspace{-.5mm}[b, \b_{_{\b_1,\b_2}}^*\hspace{-1mm}]\hspace{-1mm}+\hspace{-1mm}2tr(s\mathfrak{I}_0)\la b,\b_{_{\b_1,\b_2}}^*\hspace{-.5mm}\ra),\vspace{1mm}\\
\;[\la \b_1,\b_2\ra,v\ot
c]=\frac{1}{2\ell}\mathfrak{I}_0v\ot \b_{_{\b_1,\b_2}}^*\cdot c-\frac{1}{2}v\ot
(f(c,\b^*_2)\cdot \b^*_1+f(c,\b^*_1)\cdot \b^*_2)\\
\;[\la\b_1,\b_2\ra,\la\b'_1,\b'_2\ra]=\la d^\ell_{\b_1,\b_2}(\b'_1),\b'_2\ra+\la\b'_1,d^\ell_{\b_1,\b_2}(\b'_2)\ra
\end{array}
 \end{equation}}
for $x,y\in\gg,$  $s,t\in\ss,$ $u,v\in\v,$
$a,a'\in\aa,$ $b,b'\in\bb,$ $c,c'\in\cc,$ $\b_1,\b_2,\b_1',\b'_2\in\fb,$ is an $R-$graded  Lie algebra with grading pair $(\gg,\hh)$ where $\hh$ is the splitting Cartan subalgebra of $\gg$ defined in Lemma \ref{type-c-alg}.

$(ii)$
If
$\LL$ is  an  $R-$graded Lie algebra with grading pair $(\fg,\fh),$ then there is a coordinate quadruple $(\fa,*,\cc,f)$ of type $BC$ and a subspace $\kk$ of $\fb:=\fb(\fa,*,\cc,f)$ satisfying the uniform property on $\fb$ such that $\LL$ is isomorphic to $\LL(\fb,\kk).$

\end{Theorem}

\pf $(i)$ We prove that $\LL(\fb,\kk)$ together with
(\ref{probc-gen}) is  a Lie algebra. For $\lam\in \Lam$ set $n_\lam:=|I_\lam|$ and 
$\LL^\lam:=(\gg^\lam\dot\ot\aa)\op(\ss^\lam\dot\ot
\bb)\op(\v^\lam\dot\ot\cc)\op\la\fb,\fb\ra.$ Also for
$a,a'\in\aa,b,b'\in\bb,$ and $c,c'\in\cc,$ set
$$\begin{array}{l}\la
a,a'\ra_\lam: =(((\frac{-1}{\ell}\mathfrak{I}_0+\frac{1}{n_\lam}\mathfrak{I}_\lam)\ot\frac{1}{2}[a,a'])+\la a,a'\ra,\\
\la
b,b'\ra_\lam :=((\frac{-1}{\ell}\mathfrak{I}_0+\frac{1}{n_\lam}\mathfrak{I}_\lam)\ot\frac{1}{2}[b,b'])+\la b,b'\ra,\\
\la
c,c'\ra_\lam :=((\frac{1}{\ell}\mathfrak{I}_0-\frac{1}{n_\lam}\mathfrak{I}_\lam)\ot\frac{1}{2} c\heart c')+\la c,c'\ra,\\
\la a,b\ra_\lam=\la b,c\ra_\lam=\la a,c\ra_\lam:=0.\end{array}$$
Take $\la\fb,\fb\ra_\lam:=\hbox{span}\{\la a,a'\ra_\lam,\la b,b'\ra_\lam,\la c,c'\ra_\lam\mid a,a'\in\aa,b,b'\in\bb,c,c'\in\cc\}$ and note that  as $\kk$ satisfies the uniform property on $\fb,$ we have  $$\LL^\lam:=(\gg^\lam\dot\ot\aa)\op(\ss^\lam\dot\ot
\bb)\op(\v^\lam\dot\ot\cc)\op\la\fb,\fb\ra_\lam.$$ For $\b=a+b+c,\b'=a'+b'+c'\in\fb,$ set $\la \b,\b'\ra_\lam:=\la a,a'\ra_\lam+\la b,b'\ra_\lam+\la c,c'\ra_\lam.$

Now consider the linear transformation $\psi:\fb\ot\fb\longrightarrow\la\fb,\fb\ra_\lam $ mapping $\b\ot\b'$ to $\la\b,\b'\ra_\lam.$
It is not difficult to see that  for the  subspace $K$ of $\fb\ot\fb$ defined in Subsection \ref{subsect2-1}, $\psi(K)=\{0\}.$ So $\psi$ induces a linear transformation   $\dot\psi:\{\fb,\fb\}_{n_\lam}\longrightarrow\la\fb,\fb\ra_\lam $ mapping $\{\b,\b'\}_{n_\lam}$ to $\la\b,\b'\ra_\lam.$ Take $\kk_\lam$ to be the kernel of $\dot\psi.$ If $t\in\bbbn,$ $a_i,a'_i\in\aa, $ $b_i,b'_i\in\bb$ and $c_i,c'_i\in\cc$ ($1\leq i\leq t$) are such that $\sum_{i=1}^t(\{a_i,a'_i\}_{n_\lam}+\{b_i,b'_i\}_{n_\lam}+\{c_i,c'_i\}_{n_\lam})\in\kk_{n_\lam},$ then  $\sum_{i=1}^t(\la a_i,a'_i\ra_\lam+\la b_i,b'_i\ra_\lam+\la c_i,c'_i\ra_\lam)=0.$ This implies that  {\small \begin{equation}\label{the-final}((\frac{-1}{\ell}\mathfrak{I}_0+\frac{1}{n_\lam}\mathfrak{I}_\lam)\ot\frac{1}{2}\sum_{i=1}^t([ a_i,a'_i]+[b_i,b'_i]-(c_i\heart c'_i))+\sum_{i=1}^t(\la a_i,a'_i\ra+\la b_i,b'_i\ra+\la c_i,c'_i\ra)=0.\end{equation}} This in turn implies that $\sum_{i=1}^t(\la a_i,a'_i\ra+\la b_i,b'_i\ra+\la c_i,c'_i\ra)=0.$
Therefore we get that  $\kk_\lam$ is a subset of the full skew-dihedral homology group of $\fb$ with respect to $\ell.$ But if $\lam\neq 0,$ (\ref{the-final}) implies that   $\sum_{i=1}^t([ a_i,a'_i]+[b_i,b'_i]-(c_i\heart c'_i))=0.$ Now one gets using this  together  with the fact that $\kk_\lam$ is a subset of the full skew-dihedral homology group of $\fb$ with respect to $\ell,$ that $\kk_\lam$ is a subset of the full skew-dihedral homology group of $\fb$ with respect to $n_\lam.$ Now  it follows from \cite[Chapter $V$]{ABG2} that $\LL^\lam$  together with the  product  introduced in (\ref{probc-gen}) restricted to $\LL^\lam\times\LL^\lam$  defines a Lie algebra. Therefore  $\LL$ together with $[\cdot,\cdot]$ is a Lie algebra as $\LL=\cup_{\lam\in\Lam}\LL^\lam.$ Now one can  easily see  that $\LL$ has a weight space decomposition $\LL=\op_{\a\in R}\LL_\a$ with respect to $\hh$ in which $$\LL_\a=\left\{\begin{array}{ll}\v_\a\ot \cc& \hbox{if } \a\in R_{sh}\\
(\gg_\a\ot\aa)\op(\ss_\a\ot \bb)& \hbox{if } \a\in R_{lg}\\
\gg_\a\ot\aa&\hbox{if } \a\in R_{ex}\\
(\gg_0\ot \aa)\op(\ss_0\ot \bb)\op\la\fb,\fb\ra& \hbox{if } \a=0\end{array}\right.$$ and  that $\LL$ is an $R-$graded Lie algebra with grading pair $(\gg,\hh).$

$(ii)$
For $\lam\in\Lam,$ set
$$\begin{array}{l}\LL^\lam:=\sum_{\a\in
R_\lam^\times}\LL_\a\op\sum_{\a\in
R_\lam^\times}[\LL_\a,\LL_{-\a}],\\\\
\fg^{\lam}:=\sum_{\a\in
(R_\lam)_{sdiv}^\times}\fg_\a\op\sum_{\a\in
(R_\lam)_{sdiv}^\times}[\fg_\a,\fg_{-\a}]
\end{array}$$ and note that    $\fg^\lam$ is isomorphic to $\gg^\lam.$ We know  by Lemma \ref{final1} that $\LL^\lam$ is  an $R_\lam-$graded  Lie algebra with grading pair $(\fg^\lam,\fh^\lam:=\fg^\lam\cap\fh).$
Consider $\LL^0$ as  a $\fg^0-$module and suppose that  $\{\gg_i^0,\ss^0_j,\v^0_t,\dd_0\mid i\in \ii,j\in\jj,t\in\T\}$ is a class of finite dimensional   $\fg^0-$submodules of $\LL^0$ such that
\begin{itemize}
\item $\LL^0=\sum_{i\in \ii}\gg_i^0\op\sum_{j\in\jj}\ss^0_j\op\sum_{i\in
\T}\v^0_t\op\dd_0,$
\item $\dd_0$  is a trivial $\fg^0-$submodule of $\LL^0,$
\item  for   $i\in \ii,j\in \jj$ and $t\in \T,$ $\gg_i^0$ is isomorphic to $\gg^0,$  $\ss^0_j$ is isomorphic to $\ss^0,$ and $\v^0_t$ is isomorphic to $\v^0.$
\end{itemize}

 Now  for $\lam\in\Lam,$ consider $\LL^\lam$ as a $\fg^\lam-$module via the adjoint representation. Using Lemmas \ref{divide1} and \ref{divide2}, one finds finite dimensional irreducible $\fg^\lam-$submodules $\gg_i^\lam,$ $\ss^\lam_j,$ $\v^\lam_t$ ($i\in\ii,j\in\jj,t\in \T$) of $\LL^\lam$ and a trivial $\fg^\lam-$submodule $\dd_\lam$ such that
%
%
$$(\ii,\jj,\T,\{\gg^0_i\},\{\gg^\lam_i\},\{\ss^0_j\},\{\ss^\lam_j\},\{\v^0_t\},\{\v^\lam_t\}, \dd_0,\dd_\lam)$$ is an $(R^\lam,R^0)-$datum  for the pair $(\LL^\lam,\LL^0)$ (see (\ref{last6})).
we know from Subsection \ref{subsect2-2} that there is a coordinate quadruple $(\fa,*,\cc,f)$ of type $BC$ and a subspace $\kk_\lam$ of the full skew-dihedral homology group of $\fb:=\fb(\fa,*,\cc,f)$ with respect to $n_\lam=|I_\lam|$ such that   $\dd_\lam$ is a subalgebra of $\LL^\lam$ isomorphic to  the quotient algebra $\{\fb,\fb\}_{n_\lam}/\kk_\lam,$ say via    $\phi_\lam:\{\fb,\fb\}_{n_\lam}/\kk_\lam\longrightarrow\dd_\lam.$  Now for $\b,\b'\in\fb,$ set \begin{equation}\label{2setare}\la\b,\b'\ra^\lam:=\phi_\lam(\{\b,\b'\}_\lam+\kk_\lam).\end{equation} Take $\aa$ and $\bb$ to be the $*-$fixed and $*-$skew fixed points of $\fa$ respectively and note that $$\dd_\lam=\hbox{span}\{\la a,a'\ra^\lam,\la b,b'\ra^\lam,\la c,c'\ra^\lam\mid a,a'\in\aa,b,b'\in\bb,c,c'\in\cc\}.$$
We now proceed with  the proof in the following  steps:

\medskip

\noindent\underline{\textbf{Step 1:}} For  $i\in \ii,$ $j\in \jj,$ $t\in \T$ and $\lam,\mu\in\Lam$ with $\lam\preccurlyeq\mu,$
$\gg_i^\mu$ is the $\fg^\mu-$submodule of $\LL^\mu$ generated by $\gg_i^\lam,$ $\ss_j^\mu$ is the $\fg^\mu-$submodule of $\LL^\mu$ generated by  $\ss_j^\lam$ and
$\v^\mu_t$ is the $\fg^\mu-$submodule of $\LL^\mu$ generated by $\v_t^\lam.$ In other words, $$(\ii,\jj,\T,\{\gg^\lam_i\},\{\gg^\mu_i\},\{\ss^\lam_j\},\{\ss^\mu_j\},\{\v^\lam_t\},\{\v^\mu_t\}, \dd_\lam,\dd_\mu)$$ is an $(R^\mu,R^\lam)-$datum for the pair $(\LL^\mu,\LL^\lam):$ It is immediate using the facts that $\fg^\lam$  is a subalgebra of $\fg^\mu.$
%
\bigskip

\noindent\underline{\textbf{Step 2:}} For $\lam\in\Lam,$
$\LL^\lam=\sum_{i\in \ii}\gg_i^\lam\op\sum_{j\in\jj}\ss^\lam_j\op\sum_{i\in
\T}\v^\lam_t\op\dd_0:$   By Step 1 and  Remark \ref{rem1},   $\LL^\lam=(\sum_{i\in \ii}\gg_i^\lam\op\sum_{j\in\jj}\ss^\lam_j\op\sum_{i\in
\T}\v^\lam_t)+\dd_0.$ Suppose $d\in\dd_0,$  $x\in\sum_{i\in \ii}\gg_i^\lam\op\sum_{j\in\jj}\ss^\lam_j\op\sum_{i\in
\T}\v^\lam_t$ and  $x+d=0.$
Since $d\in \dd_0,$  there are $t\in\bbbn,$ $a_i,a'_i\in\aa,$ $b_i,b'_i\in\bb$ and $c_i,c'_i\in\cc$  $(1\leq i\leq t)$ such that $d=\sum_{i=1}^t\la a_i,a'_i\ra^0+\la b_i,b'_i\ra^0+\la c_i,c'_i\ra^0.$
It follows from  Step 1 and Lemma \ref{divide4} that  there is $y\in\sum_{i\in \ii}\gg_i^\lam\op\sum_{j\in\jj}\ss^\lam_j\op\sum_{i\in
\T}\v^\lam_t$ such that $d=y+\sum_{i=1}^t\la a_i,a'_i\ra^\lam+\la b_i,b'_i\ra^\lam+\la c_i,c'_i\ra^\lam.$ Now as $0=x+d=x+y+\sum_{i=1}^t\la a_i,a'_i\ra^\lam+\la b_i,b'_i\ra^\lam+\la c_i,c'_i\ra^\lam,$ we get that $x+y=0$ and $\sum_{i=1}^t(\la a_i,a'_i\ra^\lam+\la b_i,b'_i\ra^\lam+\la c_i,c'_i\ra^\lam)=0.$  Take $\mu\in \Lam $ to be such that $\lam\preccurlyeq\mu,$ then using Step 1, one gets that the  pairs $(\LL^\lam,\LL^\mu)$  and $(\LL^0,\LL^\lam)$ play the same role as  the pair $(\LL^\ell,\LL^n)$ in Subsection \ref{subsect2-2}. Using   Remark \ref{rem1} for the pair $(\LL^\lam,\LL^\mu),$ one gets that   $\sum_{i=1}^t(\la a_i,a'_i\ra^\mu+\la b_i,b'_i\ra^\mu+\la c_i,c'_i\ra^\mu)=0$ and $\sum_{i=1}^t([a_i,a'_i]+[b_i,b'_i]-c_i\heart c'_i)=0.$ Next using Remark \ref{rem1} for  the pair $(\LL^0,\LL^\mu),$ we get that $d=\sum_{i=1}^t(\la a_i,a'_i\ra^0+\la b_i,b'_i\ra^0+\la c_i,c'_i\ra^0)=0.$
\bigskip

\noindent\underline{\textbf{Step 3:}} $\kk_0$ satisfies the uniform property  on $\fb:$ Suppose that $$\sum_{i=1}^t(\{a_i,a'_i\}_\ell+\{b_i,b'_i\}_\ell+\{c_i,c'_i\}_\ell)\in\kk_0,$$ for  $a_1,a'_1,\ldots,a_t,a'_t\in\aa,$ $b_1,b'_1,\ldots,b_t,b'_t\in\bb,$ and  $c_1,c'_1,\ldots,c_n,c'_n\in\cc,$   so $\sum_{i=1}^n(\la a_i,a'_i\ra^0+\la b_i,b'_i\ra^0+\la c_i,c'_i\ra^0)=0.$ Now take $\lam\in\Lam\setminus\{0\},$ then by Step 1, $(\LL^0,\LL^\lam)$ plays the same role as  the pair $(\LL^\ell,\LL^n)$ in Subsection \ref{subsect2-2} and so an argument analogous to the proof of  Step  2 shows that $\sum_{i=1}^t([a_i,a'_i]+[b_i,b'_i]-c_i\heart c'_i)=0.$ This completes the proof.

\bigskip

\noindent\underline{\textbf{Step 4:}}
{\small$$\displaystyle{ \bigcup_{\lam\in\Lam}\sum_{i\in
\ii}\gg_i^\lam=\sum_{i\in\ii}\bigcup_{\lam\in\Lam}\gg_i^\lam},
\displaystyle{\bigcup_{\lam\in\Lam}\sum_{j\in\jj}\ss^\lam_j=\sum_{j\in\jj}\bigcup_{\lam\in\Lam}\ss^\lam_j,}\; \displaystyle{\bigcup_{\lam\in\Lam}\sum_{t\in \T}\v^\lam_t}
\displaystyle{=\sum_{t\in \T}\bigcup_{\lam\in\Lam}\v^\lam_t}:$$}

We just  prove the first equality; the other ones may be   proved using an analogous argument. Suppose that $\displaystyle{x\in \sum_{i\in\ii}\bigcup_{\lam\in\Lam}\gg_i^\lam,}$ then there are  $i_1,\ldots,i_n\in I,$ $\mu_1,\ldots,\mu_n\in \Lam$ and
$x_1\in\gg_{i_1}^{\mu_1},\ldots,x_n\in\gg_{i_n}^{\mu_n}$ with $x=x_1+\cdots+x_n.$ Take
$\mu\in\Lam$ to be such that
$\mu_t\preccurlyeq\mu$ for all $1\leq t\leq n.$ Then by Step 1, for  $1\leq t\leq n,$
$\gg^{\mu_t}_{i_t}\sub\gg_{i_t}^\mu.$ Therefore
$x=x_1+\cdots+x_n\in\sum_{i\in I}\gg^\mu_i\sub\bigcup_{\lam\in\Lam}\sum_{i\in I}\gg^\lam_i.$
This completes the proof.
\bigskip

\noindent\underline{\textbf{Step 5:}} We have
$$\bigcup_{\lam\in\Lam}\LL^\lam=\bigcup_{\lam\in\Lam}(\sum_{i\in
\ii}\gg_i^\lam)+\bigcup_{\lam\in\Lam}(\sum_{j\in\jj}\ss^\lam_j)+\bigcup_{\lam\in\Lam}\sum_{t\in
\T}(\v^\lam_t)+\dd_0:$$

Using Step 2, one gets that the left hand side of the above equality is a subset of the right hand side. So we need to show the other side inequality. Suppose $a$ is an element of the  right hand side, then there are $\mu_1,\mu_2,\mu_3,$ $x\in\sum_{i\in
\ii}\gg_i^{\mu_1}\sub\LL^{\mu_1},y\in\sum_{j\in\jj}\ss^{\mu_2}_j\sub\LL^{\mu_2},z\in\sum_{t\in
\T}\v^{\mu_3}_t\sub\LL^{\mu_3}$ and $d\in\dd_{0}\sub
\LL^{0}$ such that $a=x+y+z+d.$ Fix an upper bound $\mu\in\Lam$ for $\{0,\mu_1,\mu_2,\mu_3\},$ then
we have
$a=x+y+z+d\in\LL^\mu\sub\cup_{\lam\in\Lam}\LL^\lam.$
\bigskip

\noindent\underline{\textbf{Step 6:}} $(\sum_{i\in
\ii}\bigcup_{\lam\in\Lam}\gg_i^\lam)+(\sum_{j\in\jj}\bigcup_{\lam\in\Lam}\ss^\lam_j)+(\sum_{t\in
\T}\bigcup_{\lam\in\Lam}\v^\lam_t) +\dd_0$ is a direct sum:
Suppose that $\lam_1,\ldots,\lam_n,\mu_1,\ldots,\mu_m,\eta_1,\ldots,\eta_p\in\Lam$ and  $i_1,\ldots,i_n\in \ii,$
$j_1,\ldots,j_m\in\jj$ and  $t_1,\ldots,t_p\in\T$ are distinct. Let $x_{i_1}\in\gg^{\lam_1}_{i_1},\ldots,x_{i_n}\in\gg^{\lam_n}_{i_n},$
$y_{j_1}\in\ss_{j_1}^{\mu_1},\ldots,y_{j_m}\in\ss_{j_m}^{\mu_m},$
$z_{t_1}\in\v_{t_1}^{\eta_1},\ldots,z_{t_p}\in\v_{t_p}^{\eta_p}$ and $d\in\dd_0$ be such that  $$x_{i_1}+\cdots+x_{i_n}+y_{j_1}+\cdots+y_{j_m}+z_{t_1}+\cdots+z_{t_p}+ d=0.$$
Now take  $\lam\in \Lam$ to be an upper bound for
$\{\lam_1,\ldots,\lam_n,\mu_1,\ldots,\mu_m,\eta_1,\ldots,\eta_p\}$
with respect to the partial ordering on $\Lam,$ then we get using Steps 1,2 that
$x_{i_1}+\cdots+x_{i_n}+y_{j_1}+\cdots+y_{j_m}+z_{t_1}+\cdots+z_{t_p}+d$
is a summation in $(\op_{i\in \ii}\gg^\lam_i)\bigoplus(\op_{j\in
\jj}\ss^\lam_j)\bigoplus(\op_{t\in \T}\v^\lam_t)\op\dd_0.$ Therefore we have
$$x_{i_1}=0,\ldots,x_{i_n}=0,y_{j_1}=0,\ldots,y_{j_m}=0,z_{t_1}=0,\ldots,z_{t_p}=0,d=0.$$This completes the proof of this step.

\medskip

\noindent\underline{\textbf{Step 7:}} The assertion stated in  $(ii)$ is true: Take $\aa$ to be a vector space with a basis $\{a_i\mid i\in\ii\},$ $\bb$ to be a vector space with a basis $\{b_j\mid j\in\jj\},$ and $\cc$ to be a vector space with a basis $\{c_t\mid t\in\T\}.$
Using Steps  1-2,4-6, we get that  \begin{eqnarray*}\LL=\bigcup_{\lam\in\Lam}\LL_\lam&=&\bigcup_{\lam\in\Lam}(\sum_{i\in
\ii}\gg_i^\lam)+\bigcup_{\lam\in\Lam}(\sum_{j\in\jj}\ss^\lam_j)+\bigcup_{\lam\in\Lam}(\sum_{t\in
\T}\v^\lam_t)+\dd_0\\
&=&(\bigoplus_{i\in
\ii}\bigcup_{\lam\in\Lam}\gg_i^\lam)\op(\bigoplus_{j\in\jj}\bigcup_{\lam\in\Lam}\ss^\lam_j)\op(\bigoplus_{t\in
\T}\bigcup_{\lam\in\Lam}\v^\lam_t)\op\dd_0.
\end{eqnarray*}

Now consider $\LL$ as a $\fg-$module via the adjoint representation and for $i\in \ii,$ $j\in\jj$ and $t\in\T,$ set
$$\gg^{(i)}:=\bigcup_{\lam\in\Lam}\gg_i^\lam,\; \ss^{(j)}:=\bigcup_{\lam\in\Lam}\ss^\lam_j,\;
\v^{(t)}:=\bigcup_{\lam\in\Lam}\v^\lam_t,$$then by Propositions \ref{dir-lim-mod} and \ref{rep-local},
$\gg^{(i)}$ is a $\fg-$submodule of $\LL$ isomorphic to $\fg\simeq\gg,$ $\ss^{(j)}$ is
a $\fg-$submodule isomorphic to $\ss$ and $\v^{(t)}$ is a
$\fg-$submodule isomorphic to $\v.$ Therefore as a vector space, we
can identify $\LL$ with
$$(\gg\ot\aa)\op(\ss\ot\bb)\op(\v\ot\cc)\op\dd_0.$$
%
%
%
%
%
such that   for each $\lam\in\Lam,$ $\LL^\lam$ is identified with $$(\gg^\lam\dot\ot\aa)\op(\ss^\lam\dot\ot\bb)\op(\v^\lam\dot\ot\cc)\op\dd_0.$$
Now for $\lam\in\Lam,$ $(\LL^\lam,\LL^0)$ plays the same role as $(\LL^n,\LL^\ell)$ in \S \ref{subsub1} and so  we are done using Step 3 together with Proposition  \ref{divide5}.\qed
\begin{Theorem}\label{type-a-d}
Suppose that $I$ is an infinite index set and $I_0$ is a subset of $I$ of cardinality $\ell>5.$ Let  $R$ be an irreducible locally finite root system of type $X=D_I$ or $X=\dot{A}_I.$  Suppose that $\v$ is a vector space with a basis $\{v_i\mid i\in I\}$ and take $\gg$ to be the  finite dimensional split simple Lie algebra of type $X$ as in Lemmas \ref{type-a-alg} or  \ref{type-d-alg} respectively. Also define $$\mathfrak{I}_0:\v\longrightarrow\v\;\;\;\;
v_i\mapsto\left\{\begin{array}{ll}v_i& i\in I_0\\
0&\hbox{otherwise.}\end{array}\right.$$
For $x,y\in\gg,$ define $$x\circ y:=xy+yx-\frac{2tr(xy)}{l+1}\mathfrak{I}_0.$$

Suppose that $(\aa,id_{\aa},\{0\},{\bf 0})$ is a coordinate quadrable of type $X$  and  $\kk$  is a subset of the full skew-dihedral homology group of $\aa$   satisfying the  uniform property on $\aa.$ Set  $$\LL(\aa,\kk):=(\gg\ot\aa)\op\la\aa,\aa\ra,$$ in which $\la\aa,\aa\ra$ is the quotient space $\{\aa,\aa\}_\ell/\kk$ (see Subsection \ref{subsect2-1}) and for $a,a'\in\aa,$ take $\la a,a'\ra:=\{a,a'\}_\ell+\kk,$ then $\LL(\aa,\kk)$ together with
\begin{equation}\label{proa-d-gen}
\begin{array}{l}
\hbox{\small$[x\ot a,y\ot a']=$}\left\{\begin{array}{ll}\hbox{\small$\;[x,y]\ot\frac{1}{2}(a\circ a')+ (x\circ y)\ot\frac{1}{2}[a,a']+tr(xy)\la a,a'\ra$}& \hbox{\small $X=\dot A_I,$}\vspace{1mm}
\\\hbox{\small$\;[x,y]\ot aa'+tr(xy)\la a,a'\ra$}&\hbox{\small$X=D_I,$}\end{array}\right.
\\
\hbox{\small$[\la
a_1,a_2\ra,x\ot a]=$}\left\{\begin{array}{ll}
\hbox{\small$-\frac{1}{2(\ell+1)}((x\circ
Id_{_{\v^\ell}})\ot[a,[a_1,a_2]]$}\\
\hbox{\small $ +[x,Id_{_{\v^\ell}}]\ot (a\circ [a_1,a_2])+2tr(Id_{_{\v^\ell}}x)\la a,[a_1,a_2]\ra)$},&\hbox{\small $X=\dot A_I,$}\\
0&\hbox{\small $X=D_I,$}\end{array}\right.\vspace{1mm}\\
\hbox{\small $[\la a_1,a_2\ra,\la a'_1,a'_2\ra]=$}\left\{\begin{array}{ll}\hbox{\small $\la d^{\ell,\aa}_{a_1,a_2}(a'_1),a'_2\ra+\la a'_1,d^{\ell,\aa}_{a_1,a_2}(a'_2)\ra,$}& \hbox{\small $X=\dot A_I,$}\\
0&\hbox{\small $X=D_I,$}\end{array}\right.
\end{array}
 \end{equation}
for $x,y\in\gg,$
$a,a',a_1,a_2,a'_1,a'_2\in\aa,$
is a Lie algebra graded by $R.$ Moreover any $R-$graded Lie algebra gives rise in this manner.
\end{Theorem}
\begin{Theorem}\label{type b-c}
Suppose that $I$ is an infinite   index set and   $I_0$ is a subset of $I$ of cardinality $\ell>4.$ Take $\gg$ to be either $\mathfrak{o}_B(I)$ or $\mathfrak{sp}(I).$ Suppose that $\v$ is a vector space with a basis  $\{v_0,v_i,v_{\bar i}\mid i\in I\}$ equipped with a nondegenerate  symmetric bilinear form  $\fm$ as in (\ref{form-b-alg}) if $\gg=\mathfrak{o}_B(I)$ and it is a vector space with a basis  $\{v_i,v_{\bar i}\mid i\in I\}$      equipped with a nondegenerate skew-symmetric bilinear form  $\fm$ as in (\ref{form-c}) if $\gg:=\mathfrak{sp}(I).$ Set $$J:=\left\{\begin{array}{ll}I_0\cup\bar{I_0}\cup\{0\}& \hbox{if $\gg=\mathfrak{o}_B(I)$}\\
I_0\cup\bar{I_0}&\hbox{if $\gg=\mathfrak{sp}(I)$}\end{array}\right.$$ and define $\mathfrak{I}_0:\v\longrightarrow \v$  to be the linear transformation defined by $$v_i\mapsto \left\{\begin{array}{ll}v_i& \hbox{if $i\in J$}\\
0& \hbox{if $i\in I\cup\bar{I} \setminus J.$}\end{array}\right.
$$
Next set $\ss:=\v$ if $\gg:=\mathfrak{o}_B(I)$  and take  $\ss$ to be as in (\ref{module-s-c})  if $\gg=\mathfrak{sp}(I).$ For $e,f\in\gg\cup\ss,$ set $$e\circ f:=ef+fe-\frac{tr(ef)}{\ell}\mathfrak{I}_0.$$   Suppose that $R$ is an irreducible locally finite root system of type $X=B_I$ or $X=C_I$ and $(\fa,*,\cc,f)$
is a coordinate quadrable of type $X.$ Take $\aa$ and $\bb$ to be the set of $*-$fixed and $*-$skew fixed points of $\fa$ respectively.
For a subset  $\kk$ of the full skew-dihedral homology group of $\fa$ satisfying the  uniform property on $\fa,$ set $$\LL(\fa,\kk):=(\gg\ot\aa)\op(\ss\ot \bb)\op\la\fa,\fa\ra,$$ in which $\la\fa,\fa\ra$ is the quotient space $\{\fa,\fa\}_\ell/\kk$ (see Subsection \ref{subsect2-1}), and for $\a,\a'\in\fa,$ take $\la \a,\a'\ra:=\{\a,\a'\}+\kk.$ Then $\LL(\fa,\kk)$  together with
\begin{equation}\label{prod-gen}
\begin{array}{l}
\;[x\ot a,y\ot a']=[x,y]\ot aa'+tr(xy)\la a,a'\ra,\vspace{1mm}\\
\;[x\ot a,s\ot b]=xs\ot ab,\vspace{1mm}\\
\;[s\ot b,t\ot b']=D_{s,t}\ot f(b,b')+(s,t)\la b,b'\ra\\
\;[\la
\a_1,\a_2\ra,x\ot a]=
x\ot d^{\ell,\fa}_{\a_1,\a_2}(a),\vspace{1mm}\\
\;[\la
\a_1,\a_2\ra,s\ot b]=
s\ot d^{\ell,\fa}_{\a_1,\a_2}(b),\vspace{1mm}\\
\;[\la \a_1,\a_2\ra,\la \a'_1,\a'_2\ra]=\la d^{\ell,\fa}_{\a_1,\a_2}(\a'_1),\a'_2\ra+\la \a_1,d^{\ell,\fa}_{\a'_1,\a'_2}(\a_2)\ra.
\end{array}
 \end{equation}
(see Definition \ref{yoshii}) for $x,y\in\gg,$  $s,t\in\ss,$
$a,a'\in\aa,$ $b,b'\in\bb,$ $\a,\a',\a_1,\a_2,\a'_1,\a'_2\in\fa,$ if $\gg=\mathfrak{o}_B(I)$ and
 \begin{equation}\label{prob-c-gen}
\begin{array}{l}
\hbox{\small$[x\ot a,y\ot a']=[x,y]\ot\frac{1}{2}(a\circ a')+ (x\circ y)\ot\frac{1}{2}[a,a']+tr(xy)\la a,a'\ra,$}\vspace{1mm}\\
\hbox{\small$[x\ot a,s\ot b]=(x\circ s)\ot\frac{1}{2}[a,b]+[x,s]\ot \frac{1}{2} (a\circ b),$}  \vspace{1mm}\\
\hbox{\small$[s\ot b,t\ot b']=[s,t]\ot\frac{1}{2}b\circ b' + (s\circ t)\ot \frac{1}{2}[b,b']+tr(st)\la b,b'\ra,$} \vspace{1mm}\\
\hbox{\small$[\la \a,\a'\ra,x\ot a]=
\frac{-1}{4\ell}((x\circ
{\mathfrak{I}_0})\ot[a,\b_{\a,\a'}^*]+[x,{\mathfrak{I}_0}]\ot (a\circ \b_{\a,\a'}^*)),$}\vspace{1mm}\\
\hbox{\small$[\la \a,\a'\ra,s\ot b]=\frac{-1}{4\ell}([s,{\mathfrak{I}_0}]\hspace{-1mm}\ot\hspace{-1mm} (b\circ \b_{\a_1,\a_2}^*)\hspace{-1mm}+\hspace{-1mm}(s\circ
{\mathfrak{I}_0})\hspace{-1mm}\ot\hspace{-1mm} [b, \b_{\a_1,\a_2}^*]\hspace{-1mm}+\hspace{-1mm}2tr(s{\mathfrak{I}_0})\la b,\b^*_{\a,\a'}\ra),$}\vspace{1mm}\\
\hbox{\small$[\la \a_1,\a_2\ra,\la \a'_1,\a'_2\ra]=\la d^{\ell,\fa}_{\a_1,\a_2}(\a'_1),\a'_2\ra+\la\a'_1,d^{\ell,\fa}_{\a_1,\a_2}(\a'_2)\ra.$}
\end{array}
 \end{equation}
(see (\ref{beta*})) for $x,y\in\gg,$  $s,t\in\ss,$
$a,a'\in\aa,$ $b,b'\in\bb,$ $\a,\a',\a_1,\a_2,\a'_1,\a'_2\in\fa,$ if $\gg=\mathfrak{sp}(I)$
is a Lie algebra graded by $R.$
Moreover any $R-$graded Lie algebra gives rise in this manner.
\end{Theorem}

\vspace{2cm}

 School of Mathematics, Institute for Research in
Fundamental
Sciences (IPM), P.O. Box: 19395-5746, Tehran, Iran.\\
\\
 Department of Mathematics, University of Isfahan, Isfahan, Iran,
P.O.Box 81745-163.\\            
ma.yousofzadeh@sci.ui.ac.ir               

\end{document}